\documentclass[11pt]{article}

\usepackage[utf8]{inputenc}

\usepackage{amssymb,xcolor,xspace}
\usepackage{amsmath,graphicx,comment,bm}
\usepackage{gastex}

\usepackage{hyperref}
\usepackage{array}

\usepackage{geometry} 

\usepackage{fancyhdr} 
\usepackage{sectsty}
\allsectionsfont{\sffamily\mdseries\upshape} 

\usepackage[nottoc,notlof,notlot]{tocbibind} 
\usepackage[titles]{tocloft} 

\def\cqfd{\skip10=\parfillskip\parfillskip=0pt
\enspace\hfill\symbolecqfd\par\parfillskip=\skip10\par\medskip}
\def\symbolecqfd{\rlap{$\sqcap$}$\sqcup$}

\newtheorem{theorem}{Theorem}[section]
\newtheorem{proposition}[theorem]{Proposition}
\newtheorem{lemma}[theorem]{Lemma}
\newtheorem{corollary}[theorem]{Corollary}
\newtheorem{claim}[theorem]{Claim}

\newtheorem{pro-fact}[theorem]{Fact}

\newtheorem{pro-example}[theorem]{Example}
\newenvironment{example}{\begin{pro-example}\rm}{\cqfd\end{pro-example}}

\newtheorem{pro-remark}[theorem]{Remark}
\newenvironment{remark}{\begin{pro-remark}\rm}{\cqfd\end{pro-remark}}

\newenvironment{preuve}{\rm \trivlist \item[\hskip \labelsep{\bf
Proof.}]}{\cqfd\endtrivlist}

\def\cqfd{\skip10=\parfillskip\parfillskip=0pt
\enspace\hfill\symbolecqfd\par\parfillskip=\skip10\par\medskip}
\def\symbolecqfd{\rlap{$\sqcap$}$\sqcup$}
\def\proof{\begin{preuve}}
\def\eop{\end{preuve}}

\def\proofof#1{\rm \trivlist \item[\hskip \labelsep{\bf
Proof of~#1.}]}
\def\eopo{\cqfd\endtrivlist}

\let\phi\varphi

\DeclareMathOperator{\relab}{\textsf{norm}}

\def\inv{^{-1}}
\let\epsilon\varepsilon

\def \calC {\mathcal{C}}
\def\ocalC{\overline{\mathcal{C}}}
\def \calI {\mathcal{I}}

\def \O {\mathcal{O}}

\def\calS  {\mathcal{S}}

\def\MM{\mathbb{M}}
\def\N{\mathbb{N}}
\def \P {\mathbb{P}}
\def\PSL{\mathsf{PSL}_2(\Z)}

\def\Z{\mathbb{Z}}

\def\sk{\mathfrak{D}}
\def\cyc{\mathfrak{G}}
\def\qcyc{\textsf{q}\mathfrak{G}}
\def\rooted{\mathfrak{R}}
\def\qrooted{\textsf{q}\mathfrak{R}}

\def\probap{\mathfrak{p}}
\def\probaq{\mathfrak{q}}

\def\core{\mathsf{silh}}
\def\qcore{\mathsf{q}\textrm{-}\mathsf{silh}}

\def\path{\textsf{path}}

\def\overlap{\mathcal{P}^{\circ\!\!\circ}}
\def\nooverlap{\mathcal{P}^{\circ\circ}}

\def \calG {\mathcal{G}}

\def\ella{\ell_2}
\def\ellb{\ell_3}
\def\ka{k_2}
\def\kb{k_3}
\def\Ella{L_2}
\def\Ellb{L_3}

\def\Kb{K_3}

\title{Silhouettes and generic properties of subgroups of the modular group}

\author{
    Fr\'ed\'erique Bassino, \small{\url{bassino@lipn.fr}}\\
    \small{Universit\'e Sorbonne Paris Nord, LIPN, CNRS UMR 7030, F-93430 Villetaneuse, France}%
    \and
    Cyril Nicaud, \small{\url{cyril.nicaud@u-pem.fr}}\\
    \small{LIGM, Univ Gustave Eiffel, CNRS, ESIEE Paris, F-77454, Marne-la-Vallée, France}%
    \and
    Pascal Weil, \small{\url{pascal.weil@labri.fr}}\\
    \small{Univ. Bordeaux, CNRS, Bordeaux INP, LaBRI, UMR 5800, F-33400 Talence, France}\thanks{%
    LaBRI, Univ. Bordeaux, 351 cours de la Lib\'eration, 33400 Talence, France.}\\
    \small{CNRS, ReLaX, UMI 2000, Siruseri, India}
    }

\begin{document}

\maketitle

\begin{abstract}
We show that the probability for a finitely generated subgroup of the modular group, of size $n$, to be almost malnormal or non-parabolic, tends to 0 as $n$ tends to infinity --- where the notion of the size of a subgroup is based on a natural graph-theoretic representation of the subgroup.

The proofs of these results rely on the combinatorial and asymptotic study of a natural map, which associates with any finitely generated subgroup of $\PSL$ a graph which we call its silhouette, which can be interpreted as a conjugacy class of free finite index subgroups of $\PSL$.
\end{abstract}

\section{Introduction}

The study of the modular group $\PSL$ and of its subgroups has played a central role in algebra, number theory and geometry since the late 19th century. This paper fits in this stream of research. We point out that, while a vast literature has concentrated on the finite index subgroups of $\PSL$, we deal instead with \emph{all} finitely generated subgroups of $\PSL$.

Recall that a subgroup $H$ is \emph{almost malnormal} if $H \cap H^x$ (where $ H^x$ stands for $x\inv Hx$) is finite for every $x\not\in H$; it is \emph{parabolic} if it contains a parabolic element.
Our main results are the following: almost malnormality is negligible and parabolicity is generic for finitely generated subgroups of $\PSL$.

These results refer implicitly to a distribution of probabilities on the set of finitely generated subgroups of $\PSL$, that we now explain. With every finitely generated subgroup $H$ of $\PSL$, one associates a (computable) finite edge-labeled graph $\Gamma(H)$, called its \emph{Stallings graph}. The notion is an extension of the graphs Stallings  introduced in 1983 \cite{1983:Stallings} to represent finitely generated subgroups of free groups. The idea of using finite graphs to study subgroups of infinite groups has a long and distinguished history, with results of Gersten and Short \cite{1991:GerstenShort,1991:Short}, Arzhantseva and Ol'shanski\u{\i} \cite{1996:ArzhantsevaOlshanskii,1998:Arzhantseva}, Gitik \cite{1996:Gitik} and Kapovich \cite{1996:Kapovich} in the 1990s. The definition we use here finds its roots in Gitik's work on hyperbolic groups \cite{1996:Gitik} and in explicit constructions by Markus-Epstein \cite{2007:Markus-Epstein} for subgroups of amalgamated products of finite groups. Its exact form was given by Kharlampovich, Miasnikov and Weil \cite{2017:KharlampovichMiasnikovWeil}.

More precise statements are given in Section~\ref{sec: rappels}.
At this point, we just want to point out the following. Say that the \emph{size} of a finitely generated subgroup $H \le \PSL$ is the number of vertices in its Stallings graph\footnote{For a finite index subgroup, this \emph{size} coincides with the index of the subgroup.}. Then, for each $n \ge 1$,  $\PSL$ has only finitely many subgroups of size $n$, and we consider the uniform distribution on this finite set.\footnote{Exact and asymptotic enumeration results for the size $n$ subgroups of $\PSL$ can be found in \cite{2021:BassinoNicaudWeil,2023:BassinoNicaudWeil}.} A property of subgroups is \emph{negligible} (respectively, \emph{generic}) if the proportion of size $n$ subgroups with the property tends to 0 (respectively, 1) when $n$ tends to infinity.

Note that this randomness model strongly differs from those considered by Gilman, Miasnikov and Osin in \cite{2010:GilmanMiasnikovOsin} for subgroups of hyperbolic groups and by Maher and Sisto \cite{2019:MaherSisto} for subgroups of acylindrically hyperbolic groups. Roughly speaking, these models rely on fixing an integer $k$, randomly choosing $k$ elements of $G$ using an $n$-step Markovian mechanism, considering the subgroup generated by these $k$ elements, and letting $n$ tend to infinity. It is interesting to note that in this \emph{few generators} model,  almost malnormality is generic \cite[Theorem 1.1(2)]{2019:MaherSisto}, in contrast with our model. A similar situation is already known to arise for subgroups of free groups: in the $k$-generated model, malnormality is generic (Jitsukawa \cite{2002:Jitsukawa}), whereas in the so-called graph-based model, it is negligible (Bassino, Martino, Nicaud and Weil \cite{2013:BassinoMartinoNicaud}). Jitsukawa's result was extended also to the Gromov density model, where the number $k$ of generators is allowed to vary as an exponential function of $n$ (Bassino, Nicaud and Weil \cite{2016:BassinoNicaudWeilCM}). One can argue that the model we consider in this paper is particularly natural since Stallings graphs are in bijection with subgroups.

The proof of our results relies on a natural construction which we call the \emph{silhouetting} of the Stallings graph of a finitely generated subgroup of $\PSL$. This construction was introduced by the authors in \cite{2023:BassinoNicaudWeil}, to give an efficient random generation algorithm for subgroups of $\PSL$ of a given size and isomorphism type. The silhouette of the Stallings graph of a subgroup $H$ is, except in extremal cases, a uniform degree loop-free graph which captures essential features of the ``shape'' of $\Gamma(H)$. It is obtained by a sequence of ``simplifications'' of the graph. Our result in this paper exploits a remarkable, and somewhat surprising property of the silhouetting construction, namely the fact that it preserves uniformity. More precisely, among the size $n$ subgroups whose silhouette has size $s$, every size $s$ silhouette graph is equally likely.

This property allows us to lift asymptotic properties of silhouette graphs, which are more easily understood, to all Stallings graphs of finitely generated subgroups of $\PSL$. The proof of the uniformity preservation result is of a combinatorial nature, and relies on a fine description of the simplification operations carried out in \cite{2023:BassinoNicaudWeil}.

We would like to point out also an intermediate result which may be of independent interest. We show that, with high probability, in a finite group of permutations generated by a pair of fixpoint-free permutations $(\sigma_2,\sigma_3)$, of order 2 and 3 respectively, the composition $\sigma_2\sigma_3$ admits orbits of a certain, relatively small size (Proposition~\ref{pro:small cycles}). Obtaining such results on the composition of two randomly chosen mappings is notoriously difficult. It is, for instance, a bottleneck in the study of the properties of random deterministic automata~\cite{2014:Nicaud}. Most known results rely on a fine grain independent analysis of the mappings, but we know very little on their composition. Character theory has been used to tackle this kind of difficulties in the study of combinatorial maps~\cite{1990:JacksonVisentin,2006:Gamburd}: this approach yields enumeration results on triplets $(\sigma_1,\sigma_2,\tau)$ such that $\sigma_1\circ\sigma_2=\tau$, but for a fixed cyclic type of $\tau$ only, and it seems very difficult to exploit such results for our purposes (see also \cite{2021:BudzinskiCurienPetri}).

We note that, while there are a good number of results in the literature about the genericity or negligibility of certain properties of subgroups of free groups (\cite{2002:Jitsukawa,2013:BassinoMartinoNicaud,2008:BassinoNicaudWeil,2016:BassinoNicaudWeilCM,2016:BassinoNicaudWeil}), there are precious few such results for subgroups of other groups. We can cite in this direction, for their pioneering methods, the results of Arzhantseva and Ol'shanski\u{\i} \cite{1996:ArzhantsevaOlshanskii} who show for instance that, for a very large (generic) class of $r$-generator, $k$-relator presentations ($r,k$ fixed), all $\ell$-generated subgroups ($\ell < r$) are free and quasi-convex; and the results mentioned above of Gilman, Miasnikov and Osin \cite{2010:GilmanMiasnikovOsin} and Maher and Sisto \cite{2019:MaherSisto} on $k$-generated subgroups of a fixed hyperbolic or acylindrically hyperbolic group. To our knowledge, our results on almost malnormality  and parabolicity for subgroups of the modular group are the first that are based on the distribution of subgroups given by Stallings graphs.

\paragraph{Organization of the paper}
Section~\ref{sec: rappels} reviews the definitions of the Stallings graph of a subgroup of $\PSL$ and of the combinatorial type of a subgroup.

The silhouetting operation on Stallings graphs is introduced in
Section~\ref{sec: moves and silhouette}. More precisely, we first introduce a number of local moves in a Stallings graph in Section~\ref{sec: moves}. Iterating these moves turns out to be a confluent process, leading to the so-called \emph{silhouette} of the given graph or subgroup. We show there that silhouetting preserves uniformity (Theorems~\ref{thm:uniformity} and~\ref{thm:uniformity rooted}), and that silhouetting only moderately decreases the size of a graph: with high probability (more precisely, super-polynomially generically), the silhouette of a size $n$ subgroup has size at least $n-3n^{\frac23}$(Propositions~\ref{sec: size silhouette} and~\ref{prop: large silhouette rooted}).

The last section, Section~\ref{sec: generic properties} contains the proof of our main results: subgroups of $\PSL$ generically contain parabolic elements (Proposition~\ref{prop: parabolicity}) and fail to be almost malnormal (Theorem~\ref{thm: negligibility}). Both results exploit the statistical result on the existence of cycles reading non-trivial powers of $ab$ (where $a$ and $b$ are the order 2 and order 3 generators of $\PSL$) in the Stallings graph of a subgroup (Theorem~\ref{thm: small cycles}), whose proof reduces to proving the same (highly non-trivial) result on silhouette graphs (Proposition~\ref{pro:small cycles}). 

\medskip

To conclude this introductory section, we note that, according to the results presented here and in~\cite{2023:BassinoNicaudWeil}, the silhouetting operation is combinatorially and asymptotically significant in the study of finitely generated subgroups of $\PSL$. We are able to lift a statistical property of silhouette graphs to the class of all $\PSL$-reduced graphs, and it would be interesting to see what other properties can be lifted in that fashion. In addition, we think that the silhouetting operation also has a topological, or possibly a geometric interpretation, even for finite index subgroups, and we would be curious about its properties.

\section{Preliminaries}\label{sec: rappels}

We use the following presentation of the modular group, seen here as the free product of two cyclic groups:
$$\PSL = \Z_2 \ast \Z_3 = \langle a,b \mid a^2 = b^3 = 1\rangle.$$

The elements of $\PSL$ are represented by words over the alphabet $\{a,b,a\inv,b\inv\}$, or rather of $\{a,b,b\inv\}$ since $a\inv = a$. Each non-trivial element of $\PSL$ has a unique shortest (or \emph{normal}, or \emph{geodesic}) representative, which is a freely reduced word without factors in $\{a^2, b^2, b^{-2}\}$, that is, a word which either has length 1, or alternates letters $a$ and letters in $\{b,b\inv\}$.

An element $g\ne 1$ in $\PSL$ is said to be \emph{cyclically reduced} if it has length 1 or if its normal form starts with $a$ and ends with $b^{\pm1}$, or starts with $b^{\pm1}$ and ends with $a$. It is immediate that every non trivial element of $\PSL$ is conjugated to a cyclically reduced element. Moreover \cite[Theorem IV.2.8]{1977:LyndonSchupp}, two  conjugated cyclically reduced elements are cyclic conjugates of one another (that is: their geodesic representatives are of the form, respectively, $tt'$ and $t't$).

With each finitely generated subgroup $H$ of $\PSL$, we associate its Stallings graph as in \cite{2017:KharlampovichMiasnikovWeil,2021:BassinoNicaudWeil,2023:BassinoNicaudWeil}, which is the \emph{interesting}, or \emph{significant} part of the Schreier graph of $H$.

More precisely, first recall that the \emph{Schreier graph}, or \emph{coset graph} of $H$ is the graph with vertex set $\{Hg \mid g\in \PSL\}$, with an $a$-labeled edge from $Hg$ to $Hga$ and a $b$-labeled edge from $Hg$ to $Hgb$ for every $g\in \PSL$. For every $g\in \PSL$, we also have a $b\inv$-labeled edge from $Hgb$ to $Hg$. A \emph{path} is a finite sequence of consecutive edges, and the word spelled out by the labels of these edges is the \emph{label of the path}. In particular, if $w$ is a word on the alphabet $\{a,b,b\inv\}$ and $g\in \PSL$, $w$ labels a path from $Hg$ to $Hgw$.

It follows that a word $w$ is in $H$ if and only if it labels a \emph{cycle} at vertex $v_0 = H$ (that is: a path starting and ending at vertex $H$). The \emph{Stallings graph of $H$}, written $(\Gamma(H),v_0)$, is defined as the rooted subgraph of the Schreier graph spanned by the cycles at $v_0$ labeled by the geodesic representatives of the elements of $H$ --- that is, it consists of all the vertices and edges of the Schreier graph, which appear in these cycles.

It is immediately verified that, if $H$ has finite index in $\PSL$, say, $n$, then $\Gamma(H)$ is the whole Schreier graph of $H$, and $H$ has size $n$.

\begin{example}\label{ex: Stallings graphs}
Figure~\ref{fig: stallings graphs} shows examples of Stallings graphs. These graphs are \emph{labeled graphs}, \emph{i.e} their vertices are labeled by an initial segment of $\N$. The definition of Stallings graphs does not entail labeling vertices --- only designating a base vertex.
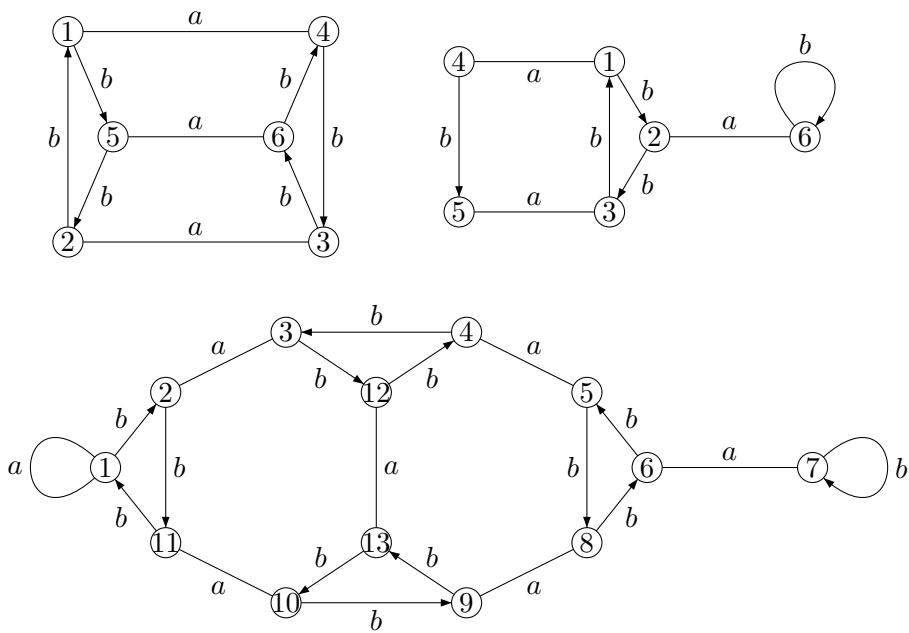
\begin{figure}[htbp]
\centering
\begin{picture}(100,82)(0,-78)
\gasset{Nw=4,Nh=4}
\node(n0)(0.0,-0.0){$1$}

\node(n1)(6.0,-14.0){5}
\node(n2)(0.0,-28.0){2}
\node(n3)(34,-0.0){4}
\node(n4)(28,-14.0){6}
\node(n5)(34.0,-28.0){3}

\drawedge(n0,n1){$b$}

\drawedge(n1,n2){$b$}

\drawedge(n2,n0){$b$}

\drawedge(n4,n3){$b$}

\drawedge(n3,n5){$b$}

\drawedge(n5,n4){$b$}

\drawedge[AHnb=0](n0,n3){$a$}

\drawedge[AHnb=0](n1,n4){$a$}

\drawedge[AHnb=0](n2,n5){$a$}
\node(m1)(72.0,-4.0){$1$}
\node(m2)(78.0,-14.0){2}
\node(m3)(72.0,-24.0){3}
\node(m4)(52,-4.0){4}
\node(m5)(52,-24.0){5}
\node(m6)(98.0,-14.0){6}

\drawedge(m1,m2){$b$}

\drawedge(m2,m3){$b$}

\drawedge(m3,m1){$b$}

\drawedge[AHnb=0](m1,m4){$a$}

\drawedge[ELside=r](m4,m5){$b$}

\drawedge[AHnb=0](m5,m3){$a$}

\drawedge[AHnb=0](m2,m6){$a$}

\drawloop(m6){$b$}
\node(h1)(5.0,-58.0){$1$}
\node(h2)(13.0,-48.0){2}
\node(h3)(29.0,-40.0){3}
\node(h4)(53.0,-40.0){4}
\node(h5)(69.0,-48.0){5}
\node(h6)(77.0,-58.0){6}
\node(h7)(99.0,-58.0){7}
\node(h8)(69.0,-68.0){8}
\node(h9)(53.0,-76.0){9}
\node(h10)(29.0,-76.0){10}
\node(h11)(13.0,-68.0){11}

\node(h12)(41.0,-48.0){12}
\node(h13)(41.0,-68.0){13}

\drawedge(h1,h2){$b$}
\drawedge(h2,h11){$b$}
\drawedge(h11,h1){$b$}

\drawedge[ELside=r](h3,h12){$b$}
\drawedge[ELside=r](h12,h4){$b$}
\drawedge[ELside=r](h4,h3){$b$}

\drawedge[ELside=r](h10,h9){$b$}
\drawedge[ELside=r](h9,h13){$b$}
\drawedge[ELside=r](h13,h10){$b$}

\drawedge[ELside=r](h5,h8){$b$}
\drawedge[ELside=r](h8,h6){$b$}
\drawedge[ELside=r](h6,h5){$b$}

\drawedge[AHnb=0](h2,h3){$a$}
\drawedge[AHnb=0](h4,h5){$a$}
\drawedge[AHnb=0](h6,h7){$a$}
\drawedge[AHnb=0](h8,h9){$a$}
\drawedge[AHnb=0](h10,h11){$a$}
\drawedge[AHnb=0](h12,h13){$a$}


\drawloop[loopangle=-180,AHnb=0](h1){$a$}
\drawloop[loopangle=0](h7){$b$}
\end{picture}
\caption{Top: The Stallings graphs of the subgroups
$H = \langle abab\inv, babab\rangle$ and $K = \langle abab, babab\inv\rangle$ of $\PSL$. Bottom: the Stallings graph of 
  $L = \langle a, (ba)^3b,  (bab\inv a)^2b, (bab\inv ab\inv a)^2b\rangle$.
In each case, the root is the vertex labeled 1.}\label{fig: stallings graphs}
\end{figure}
\end{example}

\begin{remark}\label{rk: historical}
As we noted in the introduction, finite graphs have long been used to discuss properties of arbitrary index subgroups of infinite groups. The connection with quasi-convexity and with the Howson property (which deals with question whether the intersection of finitely generated subgroups is again finitely generated) was observed by Gersten and Short in the early 1990s \cite{1991:GerstenShort,1991:Short}. Arzhantseva and Ol'shanski\u{\i} \cite{1996:ArzhantsevaOlshanskii} developed the usage of graphs to investigate the subgroups of an exponentially generic class of $k$-relator groups. The earliest usage of a graphical representation of subgroups to derive algorithmic results, may be Kapovich's work on detecting quasi-convexity in automatic groups \cite{1996:Kapovich} (note that $\PSL$ is automatic). For a systematic algorithmic approach, it is however very convenient to have a unique graphical representation, such as the Stallings graph defined above. This definition was first introduced\footnote{under the name of \emph{geodesic core}.} by Gitik \cite{1996:Gitik} for quasi-convex subgroups of hyperbolic groups (and every finitely generated subgroup of $\PSL$ is quasi-convex). It was generalized by Kharlampovich, Miasnikov and Weil \cite{2017:KharlampovichMiasnikovWeil}, who showed that the quasi-convex subgroups of automatic groups have a finite and computable Stallings graph, which led to a unified approach to a number of algorithmic decidability problems. See \cite{2017:KharlampovichMiasnikovWeil} for more historical details.
\end{remark}

\begin{remark}
The computation of the Stallings graph of a quasi-convex subgroup of an automatic group is a high complexity problem, and may not be practical in general. It is simpler when the ambient group is an amalgamated product of finite groups, as showed by Markus-Epstein \cite{2007:Markus-Epstein} (although she used slightly different graphs than our Stallings graphs). In the particular case of $\PSL$, the computation of Stallings graphs is particularly straightforward, and can be achieved in time $\O(n\log^*n)$, where $n$ is the sum of the length of the generators of the given subgroup, see \cite{2021:BassinoNicaudWeil}: given a tuple of words $\vec h = (h_1,\dots, h_r)$ on alphabet $\{a,b,b\inv\}$, one first computes the (classical) Stallings graph of the subgroup generated by $\vec h$ in the free group $F(a,b)$, as in \cite{1983:Stallings,2006:Touikan}. This is done starting with a wedge of cycles each labeled by one of the $h_i$, and then repeatedly applying Stallings foldings (in which two edges with the same label and the same start (resp. end) vertex are identified). The next step is to add an $a$-labeled edge from vertex $v$ to vertex $w$ whenever there is an $a$-labeled edge from $w$ to $v$, add a $b$-labeled edge from $v$ to $w$ whenever there is a $b^2$-labeled path from $w$ to $v$, and then apply another round of repeated Stallings foldings.
\end{remark}

It is immediate from the definition of Stallings graphs that $\Gamma(H)$ is connected and that its $a$-edges (respectively, $b$-edges) form a partial, injective map on the vertex set of the graph. Moreover, because $a^2 = b^3 = 1$, distinct $a$-edges are never adjacent to the same vertex: we distinguish therefore $a$-loops and so-called \emph{isolated $a$-edges}. Similarly, if we have two consecutive $b$-edges, say, from $v_1$ to $v_2$ and from $v_2$ to $v_3$, then $\Gamma(H)$ also has a $b$-edge from $v_3$ to $v_1$. Thus each $b$-edge is either a loop, or an \emph{isolated $b$-edge}, or a part of a $b$-triangle. Finally, every vertex except maybe the root vertex is adjacent to an $a$- and to a $b$-edge.

A rooted edge-labeled graph satisfying these conditions is called \emph{$\PSL$-reduced} and it is not difficult to see that every finite $\PSL$-reduced graph is the Stallings graph of a unique finitely generated subgroup of $\PSL$ \cite{2021:BassinoNicaudWeil}. That is, the mapping $H \mapsto (\Gamma(H),v_0)$ is a bijection between finitely generated subgroups of $\PSL$ and $\PSL$-reduced graphs. The asymptotic results that are at the heart of this paper rely on this bijection: enumeration results for subgroups are equivalent to enumeration results on $\PSL$-reduced graphs (their Stallings graphs). 

An edge-labeled graph is said to be \emph{$\PSL$-cyclically reduced} if every vertex is adjacent to an $a$- and a $b$-edge. This is equivalent to asking that $\Gamma$ be $\PSL$-reduced when rooted at every one of its vertices. We also say that a finitely generated subgroup of $\PSL$ is $\PSL$-cyclically reduced if its Stallings graph is. We note that a subgroup is $\PSL$-cyclically reduced if and only if it has minimum size in its conjugacy class (see, \textit{e.g.}, \cite[Section 2.2]{2023:BassinoNicaudWeil}.)

As is classical in combinatorics, it is actually more convenient to work with \emph{labeled} and \emph{weakly labeled graphs} (see \cite[Section II.1]{2009:FlajoletSedgewick}). A graph is said to be \emph{weakly labeled} if its vertex set is equipped with a (labeling) one-to-one map to $\N\setminus\{0\}$.\footnote{This notion of labeling, which injectively assigns an integer to each vertex, must be distinguished from the edge labeling used so far, where each edge is labeled by letter $a$ or $b$ and each path is labeled by a word.} It is \emph{labeled} if the range of the labeling map is an initial segment of $\N\setminus\{0\}$ (a set of the form $[n] = \{1,\dots,n\}$). It is immediate that a weakly labeled graph $\Gamma$ can be relabeled in a unique fashion into a labeled graph, in such a way that the order of vertices is preserved: the resulting labeled graph is called the \emph{normalization}\footnote{An operation called \emph{reduction} in \cite{2009:FlajoletSedgewick}.} of $\Gamma$, written $\relab(\Gamma)$.

To lighten up notation, we often abusively identify the vertices of a weakly labeled graph with their labels.
The graphs in Figure~\ref{fig: stallings graphs} are in fact labeled graphs.

\begin{figure}[htbp]
\centering
\begin{picture}(135,19)(3,-17)
\gasset{Nw=4,Nh=4}

\node(n0)(8.0,-5.0){}
\put(6,-15){$\Delta_1$}
\drawloop[loopangle=150,AHnb=0](n0){$a$}
\drawloop[loopangle=30](n0){$b$}

\node(n0)(33.0,-3.0){}
\node(n1)(53.0,-3.0){}
\put(41,-15){$\Delta_2$}
\drawedge[curvedepth=2](n0,n1){$b$}
\drawedge[curvedepth=2,AHnb=0](n1,n0){$a$}

\node(n0)(73.0,-5.0){}
\node(n1)(93.0,-5.0){}
\put(81,-15){$\Delta_3$}
\drawedge(n0,n1){$b$}
\drawloop[AHnb=0,loopangle=90](n0){$a$}
\drawloop[AHnb=0,loopangle=90](n1){$a$}

\node(n0)(113.0,-5.0){}
\node(n1)(133.0,-5.0){}
\put(121,-15){$\Delta_4$}
\drawedge[AHnb=0](n0,n1){$a$}
\drawloop[loopangle=90](n0){$b$}
\drawloop[loopangle=90](n1){$b$}
\end{picture}
\caption{\small All $\PSL$-cyclically reduced graphs with at most 2 vertices.}\label{fig: 2-vertex cyclically reduced}
\end{figure}
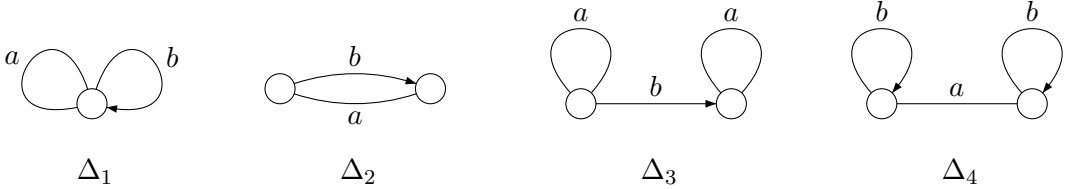

\begin{example}\label{ex: small n}
The $\PSL$-cyclically reduced graphs $\Gamma$ with 1 or 2 vertices are represented in Figure~\ref{fig: 2-vertex cyclically reduced}. 
There is only one with 1 vertex, and three with 2 vertices. Note that $\Delta_2$ and $\Delta_3$ can be labeled in two different ways while $\Delta_4$ admits only one labeling.
\end{example}

\begin{remark}\label{rk: counting convention}
Counting of graphs is done up to isomorphism. To be more precise, an isomorphism $\Gamma \to\Gamma'$ between graphs is a pair of bijections from the vertex set of $\Gamma$ to the vertex set of $\Gamma'$, and from the edge set of $\Gamma$ to the edge set of $\Gamma'$, which preserves the incidence relation. Isomorphisms between rooted graphs must also map the root of one graph to the root of the other. Finally, isomorphisms between edge-labeled graphs must also preserve these labels.
\end{remark}

If $\Gamma$ is a $\PSL$-reduced graph, its \emph{combinatorial type} is the tuple $(n,\ka,\kb,\ella,\ellb)$ where $n$ is the number of vertices of $\Gamma$, $\ka$ and $\kb$ are the numbers of isolated $a$- and $b$-edges, and $\ella$ and $\ellb$ are the numbers of $a$- and $b$-loops. We may also talk of the combinatorial type of a subgroup to mean the combinatorial type of its Stallings graph. The integer $n$ is called the \emph{size} of the graph or the subgroup. One can find  in \cite[Section 2.3.1]{2021:BassinoNicaudWeil} a discussion of the constraints on tuples that arise as combinatorial types. The combinatorial type information on a subgroup refines algebraic information such as freeness or finite index (see, e.g., \cite[Propositions 2.7, 2.9, 8.18 and Section 8.2]{2021:BassinoNicaudWeil}).

\begin{proposition}\label{prop: charact free and findex}
A subgroup $H\le \PSL$ has finite index if and only if its Stallings graph is $\PSL$-cyclically reduced and has combinatorial type of the form $(n,\ka,0,\ella,\ellb)$. It is free if and only if its combinatorial type is of the form $(n,\ka,\kb,0,0)$.

Free $\PSL$-cyclically reduced subgroups have even size, and free and finite index subgroups have size a positive multiple of $6$.
\end{proposition}

\section{Moves on $\PSL$-reduced graphs and silhouette of a subgroup}\label{sec: moves and silhouette} 

In this section, we review a combinatorial construction on labeled $\PSL$-cyclically reduced graphs first introduced in \cite{2023:BassinoNicaudWeil}. It consists in applying a number of moves on such a graph $\Gamma$, depending on its geometry, see Section~\ref{sec: moves}. These moves constitute a confluent and terminating graph rewriting system (see Section~\ref{sec: silhouettes}): this means in particular that iteratively applying these moves to a graph $\Gamma$ leads to a uniquely defined $\PSL$-cyclically reduced graph $\core(\Gamma)$, which we call the \emph{silhouette} of $\Gamma$.

This silhouetting operation was studied in \cite{2023:BassinoNicaudWeil} from a combinatorial (and, to a lesser extent, algebraic) point of view, towards exact enumeration and random generation results. We investigate in Section~\ref{sec: asymptotic properties} asymptotic properties of this operation. These properties, in turn, are used to establish certain asymptotic properties of subgroups of $\PSL$, see Section~\ref{sec: generic properties} below.

\subsection{Moves on a labeled $\PSL$-cyclically reduced graph}\label{sec: moves}
 
Except for the so-called exceptional moves defined further down, the moves we define on weakly labeled $\PSL$-cyclically reduced graphs delete vertices without changing their labels: they result in smaller weakly labeled $\PSL$-cyclically reduced graphs.

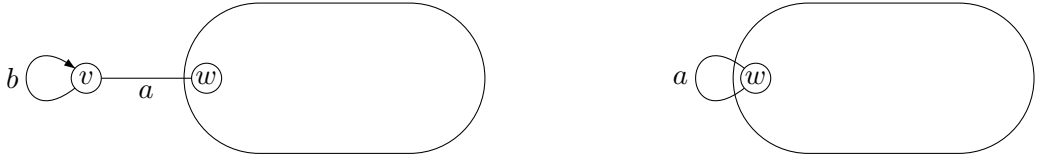
\begin{figure}[htbp]
\centering
\begin{picture}(130,20)(-2,-22)
\gasset{Nw=4,Nh=4}
\drawoval(37,-12,40,20,12)
\node(n0)(20.0,-12.0){$w$}
\node(n1)(4.0,-12.0){$v$}
\drawedge[AHnb=0](n0,n1){$a$}
\drawloop[loopangle=180,loopdiam=6](n1){$b$}

\drawoval(110,-12,40,20,12)
\node(n2)(93.0,-12.0){$w$}
\drawloop[loopangle=180,loopdiam=6,AHnb=0](n2){$a$}
\end{picture}
\caption{\small A $\lambda_3$-move}\label{fig: tb differenciation}
\end{figure}

We start with the $\lambda_3$-moves, which decrease the number of $b$-loops. Let $\Gamma$ be a weakly labeled $\PSL$-cyclically reduced graph and let $v,w$ be distinct vertices such that $\Gamma$ has a $b$-loop at $v$ and an $a$-edge between $v$ and $w$, see Figure~\ref{fig: tb differenciation}. The corresponding $\lambda_3$-move deletes vertex $v$ and the edges adjacent to it, and adds an $a$-loop at vertex $w$.

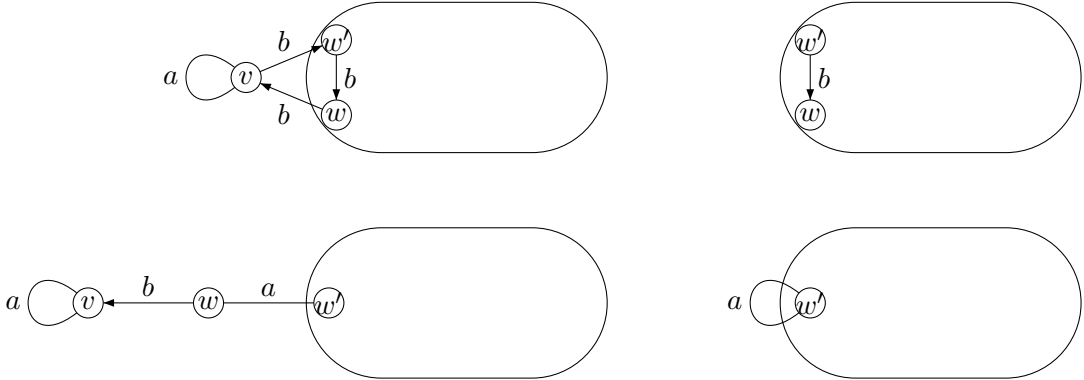
\begin{figure}[htbp]
\centering
\begin{picture}(130,50)(-2,-52)
\gasset{Nw=4,Nh=4}
\drawoval(49,-12,40,20,12)
\node(n0)(33.0,-7.0){$w'$}
\node(n00)(33.0,-17.0){$w$}
\node(n1)(21.0,-12.0){$v$}
\drawedge(n00,n1){$b$}
\drawedge(n1,n0){$b$}
\drawedge(n0,n00){$b$}
\drawloop[loopangle=180,loopdiam=6,AHnb=0](n1){$a$}

\drawoval(112,-12,40,20,12)
\node(n2)(96.0,-7.0){$w'$}
\node(n22)(96.0,-17.0){$w$}
\drawedge(n2,n22){$b$}

\drawoval(49,-42,40,20,12)
\node(n3)(32.0,-42.0){$w'$}
\node(n4)(16.0,-42.0){$w$}
\node(n5)(0.0,-42.0){$v$}
\drawedge[ELside=r](n4,n5){$b$}
\drawedge[ELside=r,AHnb=0](n3,n4){$a$}
\drawloop[loopangle=180,loopdiam=6,AHnb=0](n5){$a$}

\drawoval(112,-42,40,20,12)
\node(n6)(96.0,-42.0){$w'$}
\drawloop[loopangle=180,loopdiam=6,AHnb=0](n6){$a$}
\end{picture}
\caption{\small Above: A $\lambda_{2,1}$-move. Below: a $\lambda_{2,2}$-move}\label{fig: ta differenciation}
\end{figure}

Next we define the two kinds of $\lambda_2$-moves, which decrease the number of $a$-loops. Let $\Gamma$ be a weakly labeled $\PSL$-cyclically reduced graph and let $v$ be a vertex such that $\Gamma$ has an $a$-loop at $v$ and $v$ is part of a $b$-triangle, see Figure~\ref{fig: ta differenciation}. The corresponding $\lambda_{2,1}$-move deletes vertex $v$ and the edges adjacent to it.

The second kind of $\lambda_2$-moves is as follows. Let $\Gamma$ be a weakly labeled $\PSL$-cyclically reduced graph and let $v,w, w'$ be distinct vertices such that $\Gamma$ has an $a$-loop at $v$, an isolated $b$-edge between $v$ and $w$ (in either direction) and an $a$-edge between $w$ and $w'$, see Figure~\ref{fig: tb differenciation}. The corresponding $\lambda_{2,2}$-move deletes vertices $v$ and $w$ and the edges adjacent to them, and adds an $a$-loop at vertex $w'$.

\begin{figure}[htbp]
\centering
\begin{picture}(130,20)(-2,-22)
\gasset{Nw=4,Nh=4}
\drawoval(37,-12,40,20,12)
\node(n0)(21.0,-07.0){$v'$}
\node(n00)(21.0,-17.0){$w'$}
\node(n1)(4.0,-07.0){$v$}
\node(n11)(4.0,-17.0){$w$}
\drawedge[ELside=r](n1,n11){$b$}
\drawedge[AHnb=0](n1,n0){$a$}
\drawedge[AHnb=0](n00,n11){$a$}

\drawoval(110,-12,40,20,12)
\node(n2)(94.0,-07.0){$v'$}
\node(n22)(94.0,-17.0){$w'$}
\drawedge[AHnb=0](n2,n22){$a$}
\end{picture}
\caption{\small A $\kappa_3$-move}\label{fig: xb differenciation}
\end{figure}
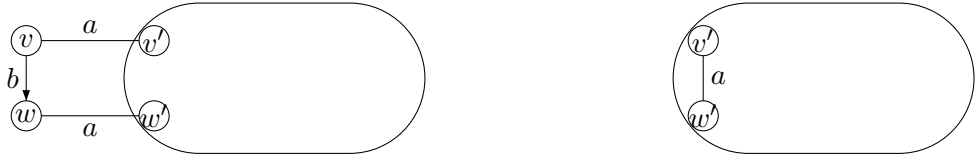

Our next kind of moves deletes isolated $b$-edges connecting vertices that sit on distinct isolated $a$-edges. More precisely, let $\Gamma$ be a weakly labeled $\PSL$-cyclically reduced graph and let $v,w,v',w'$ be distinct vertices such that $\Gamma$ has isolated $a$-edges between $v$ and $v'$, and between $w$ and $w'$, and an isolated $b$-edge between $v$ and $w$ (in either direction), see Figure~\ref{fig: xb differenciation}. The corresponding $\kappa_3$-move deletes vertices $v$ and $w$ and the edges adjacent to them, and adds an $a$-edge between $v'$ and $w'$.

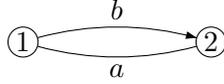
\begin{figure}[htbp]
\centering
\begin{picture}(35,9)(3,-7)
\gasset{Nw=4,Nh=4}
\node(n0)(8.0,-3.0){1}
\node(n1)(33.0,-3.0){2}
\drawedge[curvedepth=2](n0,n1){$b$}
\drawedge[curvedepth=2,AHnb=0](n1,n0){$a$}
\end{picture}
\caption{\small The preferred labeling of $\Delta_2$}\label{fig: canonical Delta2}
\end{figure}

Finally, we introduce three \emph{exceptional} moves, which can modify the labeling of vertices. One transforms an improperly labeled version of $\Delta_1$ to the (unique) labeled version of $\Delta_1$ (that is: if the unique vertex of $\Delta_1$ is not labeled by 1, the label of that vertex is made 1. Another exceptional move can be applied only to a weakly labeled version of $\Delta_3$, turning it into the labeled version of $\Delta_1$. This move can be seen as a degenerate version of a $\lambda_{2,2}$-move. The last exceptional move can be applied to any weakly labeled version of $\Delta_2$ different from the particular labeling showed in Figure~\ref{fig: canonical Delta2}, turning it to that preferred labeling of $\Delta_2$.

It is easily verified that the $\lambda_3$-, $\lambda_2$- and $\kappa_3$-moves modify the combinatorial type of a $\PSL$-cyclically reduced graph $\Gamma$ as follows. Let $\bm\tau$ be the combinatorial type of $\Gamma$.
\begin{itemize}
\item If $\Delta$ is obtained from $\Gamma$ by a $\lambda_3$-move, then $\Delta$ has combinatorial type $\bm\tau + \bm\lambda_3$, where $\bm\lambda_3 =  (-1,-1,0,1,-1)$.
\item If $\Delta$ is obtained from $\Gamma$ by a $\lambda_{2,1}$-move, then $\Delta$ has combinatorial type $\bm\tau + \bm\lambda_{2,1}$, where $\bm\lambda_{2,1} =  (-1,0,1,-1,0)$.
\item If $\Delta$ is obtained from $\Gamma$ by a $\lambda_{2,2}$-move, then $\Delta$ has combinatorial type $\bm\tau + \bm\lambda_{2,2}$, where $\bm\lambda_{2,2} =  (-2,-1,-1,0,0)$.
\item If $\Delta$ is obtained from $\Gamma$ by a $\kappa_3$-move, then $\Delta$ has combinatorial type $\bm\tau + \bm\kappa_3$, where $\bm\kappa_{3} =  (-2,-1,-1,0,0)$.
\end{itemize}

The combinatorial study of these moves carried out in \cite[Propositions 4.4, 4.5, 4.6]{2023:BassinoNicaudWeil} establishes the following enumeration results.

\begin{proposition}\label{prop: from claims}
Let $\bm\tau = (n,\ka,\kb,\ella,\ellb)$ be a combinatorial type.
\begin{itemize}
\item If $ n \ge 2$, $\ellb > 0$ and $\Delta$ is a $\PSL$-cyclically reduced graph with combinatorial type $\bm\tau + \bm\lambda_3$, then the set of labeled $\PSL$-cyclically reduced graphs $\Gamma$ of combinatorial type $\bm\tau$, such that a $\lambda_3$-move takes $\Gamma$ to $\Delta$, has $\frac{n(\ella + 1)}\ellb$ elements.

\item If $n \ge 3$, $\ella > 0$ and $\Delta$ has combinatorial type $\bm\tau + \bm\lambda_{2,1}$, then the set of labeled $\PSL$-cyclically reduced graphs $\Gamma$ of combinatorial type $\bm\tau$, such that a $\lambda_{2,1}$-move takes $\Gamma$ to $\Delta$, has $\frac{n(\kb + 1)}\ella$ elements.

\item If $n \ge 3$, $\ella > 0$ and $\Delta$ has combinatorial type $\bm\tau + \bm\lambda_{2,2}$, then the set of labeled $\PSL$-cyclically reduced graphs $\Gamma$ of combinatorial type $\bm\tau$, such that a $\lambda_{2,2}$-move takes $\Gamma$ to $\Delta$, has $2n(n-1)$ elements.

\item If $n \ge 4$, $\ella = 0$, $\kb > 0$ and $\Delta$ has combinatorial type $\bm\tau + \bm\kappa_3$, then the set of labeled $\PSL$-cyclically reduced graphs $\Gamma$ of combinatorial type $\bm\tau$, such that a $\kappa_3$-move takes $\Gamma$ to $\Delta$, has $\frac{2n(n-1)(\ka - 1)}\kb$ elements.
\end{itemize}
\end{proposition}

\subsection{Silhouette graphs, silhouette of a labeled $\PSL$-cyclically reduced graph}\label{sec: silhouettes}

Say that a $\PSL$-cyclically reduced graph $\Gamma$ (with or without a weak labeling function) is a \emph{silhouette graph} if it is equal to $\Delta_1$ or $\Delta_2$, or if it has combinatorial type $(n, n/2, 0, 0, 0)$. By Proposition~\ref{prop: charact free and findex}, the latter are the Stallings graphs of free and finite index subgroups of $\PSL$ where we forget which vertex is the base vertex, and their size $n$ is a positive multiple of 6. It is clear that, under any labeling, no move is defined on a silhouette graph of size $n > 2$.

The rewriting system on weakly labeled $\PSL$-cyclically reduced graphs, given by the moves defined in Section~\ref{sec: moves}, was showed to be \emph{confluent} in the following sense \cite[Proposition 3.4]{2023:BassinoNicaudWeil}.

\begin{proposition}\label{prop: uniqueness silhouette}
Let $\Gamma$ be a weakly labeled $\PSL$-cyclically reduced graph $\Gamma$ and let $\Delta$ and $\Delta'$ be obtained from $\Gamma$ after maximal sequences of $\lambda_3$-, $\lambda_{2,1}$-, $\lambda_{2,2}$- and $\kappa_{3}$- and exceptional moves.  Then $\Delta = \Delta'$.
\end{proposition}

In view of this result, we define the \emph{quasi-silhouette} and the \emph{silhouette} of a $\PSL$-cyclically reduced graph $\Gamma$ as $\qcore(\Gamma) = \Delta$ and $\core(\Gamma) = \relab(\Delta)$, where $\Delta$ is obtained from $\Gamma$ after a maximal sequence of moves. Proposition~\ref{prop: uniqueness silhouette} establishes that the quasi-silhouette (resp. silhouette) of $\Gamma$ is a well defined weakly labeled (resp. labeled) $\PSL$-cyclically reduced graph, which does not depend on the choice of a particular maximal sequence of moves.

\begin{example}
The first Stallings graph in Example~\ref{ex: Stallings graphs} is a silhouette graph: it is therefore its own silhouette and quasi-silhouette. The silhouette of the second Stallings graph is $\Delta_2$ with its preferred labeling (see Figure~\ref{fig: canonical Delta2}). Finally, the quasi-silhouette and silhouette of the last graph are given in Figure~\ref{fig: silhouette}.
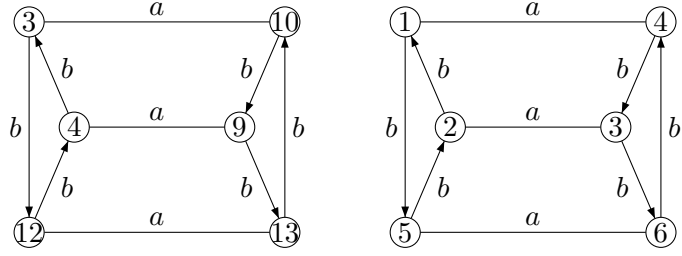
\begin{figure}[htbp]
\centering
\begin{picture}(84,34)(0,-30)
\gasset{Nw=4,Nh=4}
\node(n0)(0.0,-0.0){$3$}

\node(n1)(6.0,-14.0){4}
\node(n2)(0.0,-28.0){12}
\node(n3)(34,-0.0){10}
\node(n4)(28,-14.0){9}
\node(n5)(34.0,-28.0){13}

\drawedge[ELside=r](n0,n2){$b$}

\drawedge[ELside=r](n2,n1){$b$}

\drawedge[ELside=r](n1,n0){$b$}

\drawedge[ELside=r](n3,n4){$b$}

\drawedge[ELside=r](n4,n5){$b$}

\drawedge[ELside=r](n5,n3){$b$}

\drawedge[AHnb=0](n0,n3){$a$}

\drawedge[AHnb=0](n1,n4){$a$}

\drawedge[AHnb=0](n2,n5){$a$}
\node(m0)(50.0,-0.0){$1$}

\node(m1)(56.0,-14.0){2}
\node(m2)(50.0,-28.0){5}
\node(m3)(84,-0.0){4}
\node(m4)(78,-14.0){3}
\node(m5)(84.0,-28.0){6}

\drawedge[ELside=r](m0,m2){$b$}

\drawedge[ELside=r](m2,m1){$b$}

\drawedge[ELside=r](m1,m0){$b$}

\drawedge[ELside=r](m3,m4){$b$}

\drawedge[ELside=r](m4,m5){$b$}

\drawedge[ELside=r](m5,m3){$b$}

\drawedge[AHnb=0](m0,m3){$a$}

\drawedge[AHnb=0](m1,m4){$a$}

\drawedge[AHnb=0](m2,m5){$a$}
\end{picture}
\caption{The quasi-silhouette and the silhouette of the labeled Stallings graph of $L = \langle a, (ba)^3b,  (bab\inv a)^2b, (bab\inv ab\inv a)^2b\rangle$ given in Figure~\ref{fig: stallings graphs}.}\label{fig: silhouette}
\end{figure}
\end{example}

\begin{remark}\label{rk: order of moves}
A useful consequence of Proposition~\ref{prop: uniqueness silhouette} is the following: given a $\PSL$-cyclically reduced graph $\Gamma$, its quasi-silhouette can always be computed by iteratively performing $\lambda_3$-moves until there are no $b$-loops left, then iteratively performing $\lambda_2$-moves until there are no $a$-loops left (note that these moves do not add $b$-loops), then iteratively performing $\kappa_3$-moves until we reach a graph of size at most 2, or there are no more loops or isolated $b$-edges. In the first case, one may have to perform an exceptional move; and in the second case, a normalization step is required to obtain $\core(\Gamma)$. Several of the proofs below use this possibility of choosing the ``path'' to compute $\qcore(\Gamma)$ and $\core(\Gamma)$ in this fashion.
\end{remark}

The silhouetting operation captures algebraic properties of subgroups of $\PSL$ (see \cite[Proposition 3.4]{2023:BassinoNicaudWeil}, but we concentrate here on its probabilistic properties.

\section{Probabilistic properties of the silhouetting operation}\label{sec: asymptotic properties}

In this section, we first explore two properties of the silhouetting operation which are central to the proof of our main results. The first is that the map $\Gamma \mapsto \core(\Gamma)$ preserves uniformity on the set of labeled $\PSL$-cyclically reduced graphs of a given combinatorial type, and whose silhouette have a prescribed size, see Section~\ref{sec: randomness preserved} for the precise statement.

The second property of the silhouetting operation is that the expected number of vertices of the silhouette of a size $n$ $\PSL$-cyclically reduced graph is close to $n$, and is concentrated around this expected value, again see Section~\ref{sec: size silhouette} for a precise statement.

Finally, we show how the silhouetting operation can be extended to all finitely generated subgroups of $\PSL$ (not just the $\PSL$-cyclically reduced subgroups), and that the properties mentioned above also hold in that case.

\subsection{Silhouetting preserves uniformity on labeled $\PSL$-cyclically reduced graphs}\label{sec: randomness preserved}

If $\bm\tau = (n,\ka,\kb,\ella,\ellb)$ and $1 \le s\le n$, we let $\cyc_{\bm\tau}(s)$ (resp. $\cyc_n(s)$) be the set of labeled $\PSL$-cyclically reduced graphs of type $\bm\tau$ (resp. of size $n$) whose silhouette has size $s$. Note that $\cyc_s(s)$ is the set of labeled silhouette graphs of size $s$. Let $\P_{\bm\tau, s}$ (resp. $\P_{n,s}$) denote the uniform probability on $\cyc_{\bm\tau}(s)$ (resp. $\cyc_n(s)$). We show that the map $\core\colon \Gamma \mapsto \core(\Gamma)$, from $\cyc_{\bm\tau}(s)$ (resp. $\cyc_n(s)$) to $\cyc_s(s)$, preserves uniformity: that is, if $X$ is a subset of $\cyc_s(s)$ and $X_{\bm\tau,s}$ (resp. $X_{n,s}$) is the set of elements $\Gamma$ of $\cyc_{\bm\tau}(s)$ (resp. $\cyc_{n}(s)$) such that $\core(\Gamma)\in X$, then $\P_{\bm\tau,s}(X_{\bm\tau,s}) = \P_{n,s}(X_{n,s}) = \P_{s,s}(X)$.

\begin{theorem}\label{thm:uniformity}
Let $\bm\tau = (n,\ka,\kb,\ella,\ellb)$ be a combinatorial type and let $1\leq s\leq n$. If $\Gamma$ is an element of $\cyc_{\bm\tau}(s)$ (resp. $\cyc_n(s)$) taken uniformly at random, then $\core(\Gamma)$ is a uniformly random element of the set $\cyc_s(s)$ of size $s$ labeled silhouette graphs. That is: for any $\Delta,\Delta'\in\cyc_s(s)$, we have 
\begin{align*}
\P_{\bm\tau,s}\Big(\text{$\Gamma \in \cyc_{\bm\tau}(s)$ and }\core(\Gamma)=\Delta\Big) \enspace&=\enspace \P_{\bm\tau,s}\Big(\text{$\Gamma \in \cyc_{\bm\tau}(s)$ and }\core(\Gamma)=\Delta'\Big) \\
\P_{n,s}\Big(\text{$\Gamma \in \cyc_n(s)$ and }\core(\Gamma)=\Delta\Big) \enspace&=\enspace \P_{n,s}\Big(\text{$\Gamma \in \cyc_n(s)$ and }\core(\Gamma)=\Delta'\Big).
\end{align*}
\end{theorem}

The proof of Theorem~\ref{thm:uniformity} uses the following lemma, a ``weakly labeled version'' of the general statement. For $s > 2$, we let $\sk_n(s)$ be the set of size $s$ silhouette graphs weakly labeled with elements of $[n]$ (in particular, $\sk_s(s) = \cyc_s(s)$). If $\Delta \in \sk_n(s)$, we let $\qcyc_{\bm\tau}(\Delta)$ (resp. $\qcyc_n(\Delta)$) be the set of labeled graphs in $\cyc_{\bm\tau}(s)$ (resp. $\cyc_n(s)$) whose quasi-silhouette is $\Delta$.

\begin{lemma}\label{lm: quasi-labeled uniformity}
Let $\bm\tau = (n,\ka,\kb,\ella,\ellb)$ be a combinatorial type and let $3\leq s\leq n$. If $\Delta, \Delta' \in \sk_n(s)$, then the sets $\qcyc_{\bm\tau}(\Delta)$ and $\qcyc_{\bm\tau}(\Delta')$ have the same cardinality.
\end{lemma}

\proof
We show by induction on $n$ that the cardinality of $\qcyc_{\bm\tau}(\Delta)$ depends only on $\bm\tau$, not on $\Delta$.

If $n = s$, then $\bm\tau = (s,s/2,0,0,0)$ and $\qcyc_{\bm\tau}(\Delta) = \{\Delta\}$ has cardinality 1.

Now suppose that $n > s$. We first deal with the case where $\ellb > 0$. Let $\Gamma \in \qcyc_{\bm\tau}(\Delta)$ and let $\Gamma'$ be obtained from $\Gamma$ by a $\lambda_3$-move. Then $\qcore(\Gamma') = \qcore(\Gamma) = \Delta$ by Proposition~\ref{prop: uniqueness silhouette} and $\Gamma' \in \qcyc_{\bm\tau+\bm\lambda_3}(\Delta)$. Proposition~\ref{prop: from claims} then shows that $|\qcyc_{\bm\tau}(\Delta)| = \frac{n(\ella+1)}{\ellb}\ |\qcyc_{\bm\tau+\bm\lambda_3}(\Delta)|$, which only depends on $\bm\tau$. 

If $\ellb = 0$ and $\ella > 0$, let $v$ be a vertex carrying an $a$-loop in $\Gamma$ and let $\Gamma'$ be the $\PSL$-cyclically reduced graph obtained from $\Gamma$ by the corresponding $\lambda_2$-move: a $\lambda_{2,1}$-move if $v$ sits on a $b$-triangle of $\Gamma$, a $\lambda_{2,2}$-move otherwise. Again, $\qcore(\Gamma') = \qcore(\Gamma) = \Delta$, so $\Gamma' \in \qcyc_{\bm\tau+\bm\lambda_{2,1}}(\Delta)$ if $v$ is adjacent to a $b$-triangle and $\Gamma' \in \qcyc_{\bm\tau+\bm\lambda_{2,2}}(\Delta)$ otherwise. Proposition~\ref{prop: from claims} then establishes that 
$$|\qcyc_{\bm\tau}(\Delta)| = \frac{n(\kb+1)}{\ella}\ |\qcyc_{\bm\tau+\bm\lambda_{2,1}}(\Delta)| + 2n(n-1)\ |\qcyc_{\bm\tau+\bm\lambda_{2,2}}(\Delta)|,$$
again showing that $|\qcyc_{\bm\tau}(\Delta)|$ depends only on $\bm\tau$.

If $\ella = \ellb = 0$ and $\kb > 0$, Let $\Gamma \in \qcyc_{\bm\tau}(\Delta)$ and let $\Gamma'$ be obtained from $\Gamma$ by a $\kappa_3$-move. Here too $\qcore(\Gamma') = \qcore(\Gamma) = \Delta$ and $\Gamma' \in \qcyc_{\bm\tau+\bm\kappa_3}(\Delta)$. Proposition~\ref{prop: from claims} then shows that $|\qcyc_{\bm\tau}(\Delta)| = 2\ \frac{n(n-1)(\ka - 1)}{\kb}\ |\qcyc_{\bm\tau+\bm\kappa_3}(\Delta)|$, concluding the proof by induction.
\eop

\proofof{Theorem~\ref{thm:uniformity}}
If $s\le 2$, $\cyc_s(s)$ has only one element, so the result holds trivially. We now assume that $s\ge 3$. It suffices to show that, if $\Delta \in \cyc_s(s)$, the cardinality of the set $\{\Gamma\in \cyc_{\bm\tau}(s) \mid \core(\Gamma) = \Delta\}$ depends on $\bm\tau$ only, not on $\Delta$.

If $\Gamma'$ is a weakly labeled graph with labels in $[n]$ and $\sigma$ is a permutation on $[n]$, we denote by $\sigma(\Gamma')$ the graph obtained from $\Gamma'$ by relabeling each vertex $v$ by $\sigma(v)$. Observe that $\Gamma'$ and $\sigma(\Gamma')$ have the same combinatorial type.

Moreover, if $\Gamma\in \cyc_{\bm\tau}(s)$ and $\Delta \in \cyc_s(s)$, then $\core(\Gamma) = \Delta$ if and only if $\relab(\qcore(\Gamma)) = \Delta$, if and only if $\qcore(\Gamma) = \sigma(\Delta)$ for some permutation $\sigma$ on $[n]$ which is increasing on $[s]$. The number of such permutations depends only on $n$ and $s$ (it is exactly $\frac{n!}{s!}$).

Thus, if $\bm\tau$ is a combinatorial type, then the set $\{\Gamma\in \cyc_{\bm\tau}(s) \mid \core(\Gamma) = \Delta\}$ is the disjoint union of the $\qcyc_{\bm\tau}(\sigma(\Delta))$, where $\sigma$ runs over the permutations on $[n]$ which are increasing on $[s]$. By Lemma~\ref{lm: quasi-labeled uniformity}, all these sets have the same cardinality, which does not depend on the choice of $\Delta \in \cyc_s(s)$. Therefore the cardinality
$$|\{\Gamma\in \cyc_{\bm\tau}(s) \mid \core(\Gamma) = \Delta\}| = \frac{n!}{s!}\ |\qcyc_{\bm\tau}(\Delta)|$$
depends on $\bm\tau$ and $s$, not on $\Delta$.

The corresponding result on the cardinality of $\{\Gamma\in \cyc_n(s) \mid \core(\Gamma) = \Delta\}$ follows, since this set is the disjoint union of the $\{\Gamma\in \cyc_{\bm\tau}(s) \mid \core(\Gamma) = \Delta\}$ when $\bm\tau$ runs over the combinatorial types of size $n$.
\eopo

\subsection{Size of the silhouette of a uniform random $\PSL$-cyclically reduced graph}\label{sec: size silhouette}

Let $\Gamma$ be a labeled $\PSL$-cyclically reduced graph. The computation of $\core(\Gamma)$ shows how vertices are deleted until there are no more loops or isolated $b$-edges, or until the process has led to $\Delta_1$ or $\Delta_2$. Statistically however, only a small number of vertices are deleted. Proposition~\ref{prop: large silhouette}, a large deviation result, quantifies this statement. It follows that most vertices and edges of $\Gamma$ are untouched in the reduction to $\core(\Gamma)$. This observation will play an important role in our discussion of generic properties of subgroups in Section~\ref{sec: generic properties}, especially Theorem~\ref{thm: small cycles}.
The proof of Proposition~\ref{prop: large silhouette} relies in an essential way on results in \cite{2021:BassinoNicaudWeil}.

\begin{proposition}\label{prop: large silhouette}
If $\Gamma$ is a labeled $\PSL$-cyclically reduced graph of size $n$, the average number of vertices of $\core(\Gamma)$ is greater than $n - 2 n^{\frac23} +o(n^\frac23)$.

Moreover, if $\mu> 0$, there exists $0 < \gamma < 1$ such that the probability that $\core(\Gamma)$ has less than $n - (2+\mu) n^\frac23$ vertices is $\O(\gamma^{n^\frac13})$.
\end{proposition}

\proof
Let $\Gamma$ be a labeled $\PSL$-cyclically reduced graph with combinatorial type $(n,\ka,\kb,\ella,\ellb)$. By Proposition~\ref{prop: uniqueness silhouette}, $\core(\Gamma)$ can be computed as follows: first a sequence of $\ellb$ $\lambda_3$-moves, which delete $\ellb$ vertices (and all $\ellb$ $b$-loops) and add $\ellb$ $a$-loops.

The graph obtained after these moves has combinatorial type $(n-\ellb, \ka-\ellb, \kb,\ella+\ellb,0)$, so that $\ella+\ellb$ $\lambda_2$-moves must be performed. Every $\lambda_{2,1}$-move deletes 1 vertex and 1 $a$-loop, and adds one isolated $b$-edge. Every $\lambda_{2,2}$-moves deletes 2 vertices, 1 $a$-loop and 1 isolated $b$-edge. Therefore the $\ella+\ellb$ $\lambda_2$-moves delete at most $2(\ella+\ellb)$ vertices. They also create at most $\ella+\ellb$ isolated $b$-edges.

The graph obtained after these moves now has combinatorial type $(n', \ka', \kb',0,0)$, where $n' \le n - \ella - 2\ellb$, $\ka' \le \ka-\ella -2\ellb$ and $\kb' \le \kb + \ella + \ellb$. At that point, $\kb'$ $\kappa_3$-moves must be performed. Each deletes 2 vertices. Finally, if an exceptional move must be performed, it deletes at most 1 vertex.

Therefore the total number of vertices deleted in the computation of $\core(\Gamma)$ is at most 
$$\ellb + 2(\ella+\ellb) + 2(\kb+\ella+\ellb) = 2\kb + 4\ella + 5\ellb.$$

The statement on average values is then a direct consequence of this computation and of \cite[Proposition 5.3]{2021:BassinoNicaudWeil}.  Moreover, \cite[Theorem 5.5 and Section 7.2]{2021:BassinoNicaudWeil} establish that, for each $\mu > 0$, there exists $0 < \gamma < 1$ such that
$$
\P\left(\kb \ge \left(1+\frac\mu6\right) n^\frac23\right) = \O\left(\gamma^{n^\frac23}\right),\,
\P\left(\ella \ge 2 n^\frac12\right) = \O\left(\gamma^{n^\frac12}\right),\,
\P\left(\ellb \ge 2 n^\frac13\right) = \O\left(\gamma^{n^\frac13}\right).
$$
For $n$ large enough, $2n^\frac13 < 2n^\frac12 < \frac\mu{15}n^{\frac23}$. So
$\P\left(2\kb + 4\ella + 5\ellb \ge (2+\mu) n^\frac23\right) = \O\left(\gamma^{n^\frac13}\right)$,
which concludes the proof.
\eop

\subsection{Silhouetting all finitely generated subgroups of $\PSL$}\label{sec: silhouette all subgroups}

We have, until now, discussed the silhouetting operation for labeled $\PSL$-cyclically reduced graphs. It can be naturally extended to labeled $\PSL$-reduced graphs, and the same probabilistic results holds in that extended framework.

Let $(\Gamma,v)$ be a labeled $\PSL$-reduced graph. 
Recall that $(\Gamma,v)$ is not $\PSL$-cyclically reduced if and only $v$ is not adjacent to both an $a$- and a $b$-edge, and is the only such vertex. Let $\Gamma^\textsf{o}$ be the graph obtained from $\Gamma$ by adding an $a$-loop (resp. a $b$-loop) at $v$ if $v$ is not adjacent to an $a$-edge (resp. a $b$-edge). Then $\Gamma^\textsf{o}$ is $\PSL$-cyclically reduced\footnote{The only case where two loops must added is when $\Gamma$ is the graph with one vertex and no edge (the Stallings graph of the trivial subgroup), which we do not need to consider in the rest of this section.} and $\Gamma^\textsf{o} = \Gamma$ if $\Gamma$ is $\PSL$-cyclically reduced. We define the \emph{silhouette} $\core(\Gamma,v)$ and the \emph{quasi-silhouette} $\qcore(\Gamma,v)$ of the rooted graph $(\Gamma,v)$ to be equal to $\core(\Gamma^\textsf{o})$ and $\qcore(\Gamma^\textsf{o})$.

If $1\leq s, \ell \leq n$, we let $\rooted_{n,\ell}(s)$ be the set of size $n$ labeled $\PSL$-reduced graphs $(\Gamma,v)$ whose silhouette has $s$ vertices, and such that $\Gamma^\textsf{o}$ has $\ell$ loops (a single counter for both the $a$- and the $b$-loops). Observe that 
the set $\sk(n,s)$ of quasi-silhouettes of the elements of $\cyc_n(s)$, is also the set of quasi-silhouettes of the elements of $\bigcup_\ell \rooted_{n,\ell}(s)$. A statement parallel to Theorem~\ref{thm:uniformity} holds for rooted graphs.

\begin{theorem}\label{thm:uniformity rooted}
Let $n \ge 3$ and $1\leq s, \ell \leq n$. If $(\Gamma,v)$ is an element of $\rooted_{n,\ell}(s)$ taken uniformly at random, then $\core(\Gamma,v)$ is a uniformly random element of $\cyc_s(s)$. That is: for any $\Delta,\Delta'\in\cyc_s(s)$, we have 
$$\P(\core(\Gamma,v)=\Delta)= \P(\core(\Gamma,v)=\Delta'),$$
where $\P$ denotes the uniform probability on $\rooted_{n,\ell}(s)$.
\end{theorem}

\proof
If $\Delta \in \sk_n(s)$, let $\qrooted_{n,\ell}(\Delta)$ be the set of all $(\Gamma,v) \in \qrooted_{n,\ell}(s)$ such that $\qcore(\Gamma,v) = \Delta$. This set is the disjoint union of the sets $\qrooted_{n,\ell}(\Delta,\alpha)$ ($\alpha \in \{a,b,0\}$), where $\qrooted_{n,\ell}(\Delta,0)$ consists of the elements $(\Gamma,v) \in \qrooted_{n,\ell}(\Delta)$ such that $\Gamma$ is a labeled $\PSL$-cyclically reduced graph, and $\qrooted_{n,\ell}(\Delta,a)$ (resp. $\qrooted_{n,\ell}(\Delta,b)$) consists of the $(\Gamma,v) \in \qrooted_{n,\ell}(\Delta)$ such that $v$ is not adjacent to an $a$-loop (resp. a $b$-loop).

We have $|\qrooted_{n,\ell}(\Delta,0)| = n\,\sum_{\bm\tau} |\qcyc_{\bm\tau}(\Delta)|$  and
$|\qrooted_{n,\ell}(\Delta,a)| = \sum_{\ell_i = 0}^\ell\ell_i\,\sum_{\bm\tau} |\qcyc_{\bm\tau}(\Delta)|$ for $i=2,3$.
hence
$|\qrooted_{n,\ell}(\Delta)| = (n+\ell)\,\sum_{\bm\tau} |\qcyc_{\bm\tau}(\Delta)|$,
where the summations over $\bm\tau$ are over all combinatorial types $(n,\ka,\kb,\ella,\ellb)$ such that $\ella+\ellb = \ell$.

By Lemma~\ref{lm: quasi-labeled uniformity}, $|\qcyc_{\bm\tau}(\Delta)|$ does not depend on the choice of $\Delta$ in $\sk_n(s)$, so the quantities above depend only on $n$ and $\ell$. As a result, if $\Delta, \Delta' \in \sk_n(s)$, then $\P(\qcore(\Gamma,v)=\Delta)= \P(\qcore(\Gamma,v)=\Delta')$ (over $\qrooted_{n,\ell}(s)$).

Next, as in the proof of Theorem~\ref{thm:uniformity}, we verify that this implies that, if $\Delta, \Delta' \in \cyc_s(s)$, then $\P(\core(\Gamma,v)=\Delta)= \P(\core(\Gamma,v)=\Delta')$ (over $\rooted_{n,\ell}(s)$), which concludes the proof.
\eop

We also establish a result similar to Proposition~\ref{prop: large silhouette}.

\begin{proposition}\label{prop: large silhouette rooted}
If $(\Gamma,v)$ is a random labeled $\PSL$-reduced graph of size $n$, the average number of vertices of $\core(\Gamma,v)$ is greater than $n - 2 n^\frac23 +o(n^\frac23)$.

Moreover, if $\mu> 0$, there exists $0 < \gamma < 1$ such that the probability that $\core((\Gamma,v))$ has less than $n - (2+\mu) n^\frac23$ vertices is $\O(\gamma^{n^{\frac13}})$.
\end{proposition}

\proof
Let $\cyc_n$ and $\rooted_n$ denote the set of size $n$ labeled $\PSL$-cyclically reduced graphs and $\PSL$-reduced graphs, respectively.

Let $(\Gamma,v)\in \rooted_n$ be a $\PSL$-reduced graph. We note that $\Gamma^\textsf{o}$ has the same number of vertices and isolated $a$- and $b$-edges as $\Gamma$. If $\Gamma$ has $\ella$ $a$-loops and $\ellb$ $b$-loops, then $\Gamma^\textsf{o}$ has $\ella'$ and $\ellb'$ such loops, with $(\ell'_2,\ell'_3) = (\ell_2,\ell_3)$, $(\ell_2+1,\ell_3)$ or $(\ell_2,\ell_3+1)$.

Since we need to manipulate simultaneously random variables defined on different probability spaces (namely the uniform distributions on $\cyc_n$ and $\rooted_n$), we introduce the following notation.
\begin{itemize}
\item $\P_{\cyc_n}$ and $\P_{\rooted_n}$ denote the uniform probabilities on $\cyc_n$ and $\rooted_n$, respectively.
\item $\Ella$ (respectively, $\Ellb$, $\Kb$) is the random variable defined on $\cyc_n$ such that $\Ella(\Gamma)$ (respectively, $\Ellb(\Gamma)$, $\Kb(\Gamma)$) is the number of $a$-loops (respectively, $b$-loops, isolated $b$-edges) in $\Gamma$.
\item $\Ella'$, $\Ellb'$ and $\Kb'$ are the random variables defined on $\rooted_n$ by $\Ella'(\Gamma,v) = \Ella(\Gamma^\textsf{o})$, $\Ellb'(\Gamma,v) = \Ellb(\Gamma^\textsf{o})$ and $\Kb'(\Gamma,v) = \Kb(\Gamma^\textsf{o})$.
\end{itemize}

We know from the proof of Proposition~\ref{prop: large silhouette}, that the number of vertices of $\Gamma$ deleted when silhouetting $(\Gamma,v) \in \rooted_n$, is at most
$2\Kb'(\Gamma,v) + 4 \Ella'(\Gamma,v) + 5\Ellb'(\Gamma,v).$
So we want to show that there exists $0 < \gamma < 1$ such that
$$\P_{\rooted_n}\left(2\Kb'(\Gamma,v) + 4 \Ella'(\Gamma,v) + 5\Ellb'(\Gamma,v) \ge (2+\mu) n^\frac23\right) = \O\left(\gamma^{n^\frac13}\right).$$

Let $\Gamma' \in \cyc_n$ be a size $n$ labeled $\PSL$-cyclically reduced graph. Then $\Gamma' = \Gamma^\textsf{o}$ for exactly $n + \Ella(\Gamma') + \Ellb(\Gamma')$ values of $(\Gamma,v) \in \rooted_n$. Moreover $n \le n + \Ella(\Gamma') + \Ellb(\Gamma') \le 2n$ and $|\rooted_n|\geq n |\cyc_n|$, hence
$$\P_{\rooted_n}(\Ella'(\Gamma,v)=i) \enspace=\enspace \frac{|\{(\Gamma,v)\in\rooted_n \mid \Ella'(\Gamma,v) =i\}|}{|\rooted_n|} \enspace\leq\enspace 2 \, \P_{\cyc_n}(\Ella(\Delta) = i),$$
where the first probability is in $\rooted_n$ and the second in $\cyc_n$.
Similarly we have
\[
\P_{\rooted_n}(\Ellb'(\Gamma,v)=i) \enspace\leq\enspace  2 \,\P_{\cyc_n}(\Ellb(\Delta)=i) \text{ and }
\P_{\rooted_n}(\Kb'(\Gamma,v)=i) \enspace\leq\enspace 2 \, \P_{\cyc_n}(\Kb(\Delta)=i).
\]
Summing over successive values of $i$, this yields
\begin{align*}
& \P_{\rooted_n}\left(\Ella'(\Gamma,v)\ge 2n^{\frac12}\right) \leq   2 \,\P_{\cyc_n}\left(\Ella(\Delta) \ge 2n^{\frac12}\right), \\
& \P_{\rooted_n}\left(\Ellb'(\Gamma,v) \ge 2n^{\frac12}\right) \leq  2 \,\P_{\cyc_n}\left(\Ellb(\Delta) \ge 2n^{\frac13}\right) \\
\textrm{and } &\P_{\rooted_n}\left(\Kb'(\Gamma,v) \ge \left(1+\frac\mu6\right)n^{\frac23}\right) \enspace\leq\enspace  2 \, \P_{\cyc_n}\left(\Kb(\Delta)  \ge \left(1+\frac\mu6\right)n^{\frac23}\right).
\end{align*}

We saw in the proof of Proposition~\ref{prop: large silhouette} that the probabilities in the right-hand side terms are small, and they remain small even when multiplied by 2. In particular, there exists $0 < \gamma < 1$ such that
$$ \P_{\rooted_n}\left(\Ella'(\Gamma,v)\ge 2n^{\frac12}\right) =  \O\left(\gamma^{n^\frac12}\right),\quad
\P_{\rooted_n}\left(\Ellb'(\Gamma,v) \ge 2n^{\frac12}\right) =  \O\left(\gamma^{n^\frac13}\right)$$
$$\text{ and }
\P_{\rooted_n}\left(\Kb'(\Gamma,v) \ge \left(1+\frac\mu6\right)n^{\frac23}\right) =  \O\left(\gamma^{n^\frac23}\right).$$
As in the proof of Proposition~\ref{prop: large silhouette}, this leads to
$$\P_{\rooted_n}\left(2\Kb'(\Gamma,v) + 4 \Ella'(\Gamma,v) + 5\Ellb'(\Gamma,v) \ge (2+\mu) n^\frac23\right) = \O\left(\gamma^{n^\frac13}\right),$$
which concludes the proof.
\eop

\section{Generic properties of subgroups of $\PSL$}\label{sec: generic properties}

Recall that the \emph{parabolic elements} of $\PSL$ are the conjugates of non-trivial powers of $ab$ and that a subgroup $H$ of $\PSL$ is \emph{non-parabolic} if it contains no parabolic element.

Recall also, $H$ is \emph{almost malnormal} if, for every $x \not\in H$, $H \cap H^x$ is finite. It is \emph{malnormal} if each of these intersections is trivial: malnormality coincides with almost malnormality if $H$ is torsion-free (e.g.\ a free subgroup of $\PSL$). 

We show that, generically, a finitely generated subgroup of $\PSL$ contains parabolic elements (Proposition~\ref{prop: parabolicity}) and fails to be almost malnormal (Theorem~\ref{thm: negligibility}). More precisely, we prove the following statements.

\begin{proposition}\label{prop: parabolicity}
Let $0 < \alpha < \frac16$. A random size $n$ subgroup (respectively, cyclically reduced subgroup) of $\PSL$ is non-parabolic with probability $\O(n^{-\alpha})$.
\end{proposition}

\begin{theorem}\label{thm: negligibility}
Let $0 < \alpha < \frac16$. The probability that a size $n$ subgroup (respectively, cyclically reduced subgroup) of $\PSL$ is almost malnormal is $\O(n^{-\alpha})$.
\end{theorem}

The proofs of both statements rely on Theorem~\ref{thm: small cycles} below. This technical statement deals with the presence of cycles in a $\PSL$-reduced graph labeled by a power of $ab$. In Theorem~\ref{thm: small cycles} and in the rest of the paper, a cycle labeled by $(ab)^m$ ($m \ge 1$) is called an \emph{$ab$-cycle of size $m$}. It is clear that such a cycle has length $2m$.

\begin{theorem}\label{thm: small cycles}
  Let $0 < \alpha < \frac16$. Then a random size $n$ labeled $\PSL$-reduced (respectively, $\PSL$-cyclically reduced) graph fails to have an $ab$-cycle of size at least $2$ and at most $n^\alpha$ with probability $\O(n^{-\alpha})$.
\end{theorem}

The proof of Proposition~\ref{prop: parabolicity} is a direct application of Theorem~\ref{thm: small cycles}, and is given in Section~\ref{sec: proof parabolicity}. The proof of Theorem~\ref{thm: negligibility} is less direct and relies on a graph-theoretic characterization of almost malnormality, see Section~\ref{sec: almost malnormality}. Finally, the complex proof of Theorem~\ref{thm: small cycles} is given in Section~\ref{sec: ab-cycles in general} except for a technical but deep lemma, which is interesting in its own right and is established in Section~\ref{sec: cycles in silhouette}.

\subsection{Non-parabolic subgroups of $\PSL$} \label{sec: proof parabolicity}

Recall from Section~\ref{sec: rappels}, the notion of a cyclically reduced element of $\PSL$.

\begin{lemma}\label{lemma: cyclically reduced cycles}
Let $H$ be a finitely generated subgroup of $\PSL$ and  $(\Gamma(H),v_0)$ its Stallings graph, let $g$ be a non-trivial element of $H$. If $x$ is a word with minimal length such that $u = x\inv gx$ is cyclically reduced, then $\Gamma(H)$ has an $x$-labeled path from $v_0$ to some vertex $v_1$ and a $u$-labeled cycle at $v_1$.
\end{lemma}

\proof
The result is trivial if $g$ is cyclically reduced since, in that case, $x = 1$. Let us now assume that $g$ is not cyclically reduced, so that $x \ne 1$.

If $|u| = 1$, then $xux\inv$ is the normal form of $g$. Indeed, let $c$ be the last letter of $x$, say, $x = yc$. If $u = a$ and $c = b^{\pm1}$, or $u = b^{\pm1}$ and $c = a$, then $xux\inv$ is clearly in normal form. If instead $u = c$ or $u = c\inv$, then $yuy\inv = xux\inv = g$, and $y\inv gy = u$ is cyclically reduced, contradicting the minimality of $|x|$. It follows that $\Gamma(H)$ has a cycle at $v_0$ labeled $xux\inv$ and therefore, as announced, it has a $u$-cycle at $v_1$, where $v_1$ is the extremity of the $x$-labeled path starting at $v_0$.

If $|u| > 1$, without loss of generality, set $u = au'b^\epsilon$ (for some $\epsilon  = \pm 1$), where $u'$ is either empty or ends with an $a$. If $x = 1$, then $\Gamma(H)$ has a $u$-cycle at $v_1 = v_0$. If $x \ne 1$ and $x$ ends with an $a$, say, $x = ya$, we have $g = xux\inv = yu'b^\epsilon ay\inv$ and $u'b^\epsilon a$ is cyclically reduced, contradicting the minimality of $|x|$. So $x$ cannot end with an $a$. The same reasoning shows that it cannot end with $b^\epsilon$. So $x = yb^{-\epsilon}$. The word $ xux\inv = y b^{-\epsilon}au'b^{-\epsilon} y\inv$ is therefore the normal form of $g$. Let $v_1$, $v_2$ and $v_3$ be the vertices reached from $v_0$ after reading $x = yb^{-\epsilon}$, $y$ and $yb^\epsilon$, respectively. Then $v_1$, $v_2$ and $v_3$ constitute a $b$-triangle and there is a path labeled $au'$ from $v_1$ to $v_3$, so that Considering the cycle at $v_0$ labeled by $y b^{-\epsilon}au'b^{-\epsilon} y\inv$, we see that $v$ and $v_1$ sit on the same $b$-triangle, and that $u = au'b^\epsilon$ labels a cycle at $v_1$.
\eop

\proofof{Proposition~\ref{prop: parabolicity}}
Recall that the parabolic elements of $\PSL$ are the conjugates of the non-trivial powers of $ab$. Such an element is a product of the form $g = x\inv(ab)^mx$ ($m\in \Z$, $m\ne 0$), and we observe that $(ab)^m$ is cyclically reduced. Lemma~\ref{lemma: cyclically reduced cycles} shows that, if $H$ is a subgroup of $\PSL$, then $g \in H$ if and only if $(ab)^m$ labels a cycle in $\Gamma(H)$.

And from Theorem~\ref{thm: small cycles}: for $0 < \alpha < 1/6$, $\Gamma(H)$ contains a cycle labeled $(ab)^m$ for some $2 \le m \le n^\alpha$ with probability $1 - \O(n^\alpha)$. As a result, $H$ contains a parabolic element with probability $1 - \O(n^\alpha)$.
\eopo

\subsection{Almost malnormal subgroups of $\PSL$}\label{sec: almost malnormality}

We start with the following characterization of almost malnormality, very similar to the classical characterization of malnormality for subgroups of free groups due to Jitsukawa \cite{2002:Jitsukawa}.

\begin{proposition}\label{prop: charact malnormal}
Let $H$ be a finitely generated subgroup of $\PSL$ and let $(\Gamma,v_0)$ be its Stallings graph. Then
$H$ is almost malnormal if and only if there does not exist distinct vertices $p$ and $q$ in $\Gamma(H)$ and an element $g\in \PSL$ (in normal form) of infinite order (that is: $g$ is not a conjugate of $a$ or $b$) which labels cycles in $\Gamma(H)$ at both $p$ and $q$.
\end{proposition}

\proof
Let $h\not\in H$ such that $H \cap H^h$ is infinite. It is well known that $\PSL$ is locally quasi-convex and that such groups satisfy the Howson property \cite{1991:Short}. As a result, $H \cap H^h$ is finitely generated. In particular it is a quasi-convex subgroup and hence itself a hyperbolic group. As a result, $H \cap H^h$ contains an element $g$ of infinite order. Let $w$ be its normal form. By definition, $w$ labels a cycle in $\Gamma(H)$ at $v_0$. In addition note that the finite order elements of $\PSL$ are exactly the conjugates of $a$ and $b$.

If $g = y^p$ for some reduced word $y$ and $p > 1$, let $x$ be a minimum length word such that $u = x\inv yx$ is cyclically reduced. Then $g$ has normal form $xu^px\inv$, and $u^p$ is cyclically reduced. In particular, $x$ labels a path in $\Gamma(H)$ from $v_0$ to a vertex $v_1$, and $u^p$ labels a cycle $C$ at $v_1$. It follows that $u^p$ labels a cycle at $p$ distinct vertices along $C$. We now assume that $g$ is not a proper power.

Let $x$ be a word of minimal length such that $u = x\inv gx$ is cyclically reduced. Lemma~\ref{lemma: cyclically reduced cycles} shows that $x$ labels a path in $\Gamma(H)$ from $v_0$ to a vertex $v_1$, and that $u$ labels a cycle at $v_1$.

Since $hgh\inv \in H$, Lemma~\ref{lemma: cyclically reduced cycles} shows that $\Gamma(H)$ also has a $y$-path from $v_0$ to a vertex $v_2$ and a cycle at $v_2$ labeled by a cyclically reduced word $w$ such that $hgh\inv = ywy\inv$. Moreover, since $hgh\inv$ is conjugated to $g$, the words $u$ and $w$ are cyclic conjugates, say, $u = tt'$ and $w = t't$. Let then $v_3$ be the vertex reached from $v_2$ reading $t'$. There are $tt'$-cycles at $v_1$ and $v_3$ and we need to show that $v_1 \ne v_3$.

If $v_1 = v_3$, then $tx\inv \in H$. Moreover
$$hgh\inv = xux\inv = xt'tx\inv = xt\inv tt' tx\inv = (tx\inv)\inv g (tx\inv)$$
so $tx\inv h$ and $g$ commute. By the classical characterization of commuting elements in free products \cite[Theorem 4.5]{1976:MagnusKarrassSolitar}, it follows that either $tx\inv h$ and $g$ sit in the same conjugate of $\langle a\rangle$ or $\langle b\rangle$, or $tx\inv h$ and $g$ sit in the same cyclic subgroup.

The first case is impossible since we assumed that $g$ has infinite order. Therefore there exist $z\in \PSL$ and $p,q \in \Z$ such that $g = z^p$ and $tx\inv h = z^q$. We assumed that $g$ is not a proper power, so $p = \pm 1$. Thus $tx\inv h = g^{pq} \in H$ and hence $h\in H$, a contradiction. %
\eop

The following is a very convenient corollary of Proposition~\ref{prop: charact malnormal}.

\begin{corollary}\label{cor: ab-cycles}
Let $H$ be a subgroup of $\PSL$. If there exist a word $w$ in normal form which is not a conjugate of $a$ or $b$, and an integer $m\ge 2$ such that $w^m$ labels a cycle at some vertex $p$ in the Stallings graph $\Gamma(H)$ but $w$ does not label a cycle at $p$, then $H$ fails to be almost malnormal.
\end{corollary}

\proof
The word $w^m$ labels cycles at every vertex reached from $p$ reading $w^i$, for $0 \le i < m$ and we conclude by Proposition~\ref{prop: charact malnormal}.
\eop

Our result on the probability for a subgroup to not be almost malnormal follows directly.

\proofof{Theorem~\ref{thm: negligibility}}
Theorem~\ref{thm: small cycles} shows that for $0 < \alpha < \frac16$ the Stallings graph of a size $n$ random subgroup $H$ contains a cycle labeled $(ab)^m$ (with $2 \le m \le n^\alpha$) with probability $ 1 - \O(n^{-\alpha})$. The result then follows from Corollary~\ref{cor: ab-cycles}.
\eopo

\subsection{$ab$-cycles in $\PSL$-reduced graphs}\label{sec: ab-cycles in general}

  Let $0 < \alpha < 1$. We want to show that a size $n$ $\PSL$-reduced (respectively, $\PSL$-cyclically reduced) graph has an $ab$-cycle of size $m \in [2,n^\alpha]$ with high probability. More precisely, we want to show that the probability that this is not the case, is $\O(n^{-\alpha})$.

The structure of the proof is the following: we first show that Theorem~\ref{thm: small cycles} holds when restricted to labeled silhouette graphs (this is Proposition~\ref{pro:small cycles} below); we then lift this result to all labeled $\PSL$-cyclically reduced graphs and $\PSL$-reduced graphs using the randomness preservation results of Sections~\ref{sec: randomness preserved} and~\ref{sec: silhouette all subgroups}, and the results on the expected size of the silhouette proved in Sections~\ref{sec: size silhouette} and~\ref{sec: silhouette all subgroups}.

\begin{proposition}\label{pro:small cycles}
Let $0 < \alpha < \frac16$ and let $n > 0$ be a multiple of 6. Then a random size $n$ labeled silhouette graph fails to have an $ab$-cycle of size at least $2$ and at most $n^\alpha$, with probability $\O(n^{-\alpha})$.
\end{proposition}

Besides its importance in the proof of Theorem~\ref{thm: small cycles}, Proposition~\ref{pro:small cycles} is of independent interest, as it deals with the probability of the presence of short cycles in certain permutation groups, see the discussion in the introduction. The proof of Proposition~\ref{pro:small cycles}, given in Section~\ref{sec: cycles in silhouette}, is tricky and uses different techniques than what has been used so far in this paper.

\proofof{Theorem~\ref{thm: small cycles}}
Recall that $\cyc_n$ and $\rooted_n$ denote the sets of size $n$ labeled $\PSL$-cyclically reduced graphs and $\PSL$-reduced graphs, respectively. We first deal with $\PSL$-cyclically reduced graphs.

Fix $0 < \alpha < \frac16$. Let $\probap_\alpha(n)$ be the probability that a graph in $\cyc_n$ has no $ab$-cycle of size in $[2,n^\alpha]$ and let $\probap'_\alpha(n)$ be the probability that an element of $\cyc_n$ has a silhouette of size at least $n-3n^{\frac{2}{3}}$ and has no $ab$-cycle of length in $[2,n^\alpha]$.

Proposition~\ref{prop: large silhouette} shows that there exists $0 < \gamma < 1$ such that the probability that an element of $\cyc_n$ has a silhouette with size less than $n-3n^{\frac{2}{3}}$ is $\O(\gamma^{n^{\frac{2}{3}}})$.
This yields
$$\probap_\alpha(n) \le \probap'_\alpha(n) + \O(\gamma^{n^{\frac{2}{3}}}).$$

Let now $s \ge n-3n^{\frac23}$ and $\probap'_\alpha(n,s)$ be the probability that an element of $\cyc_n(s)$ (the set of elements of $\cyc_n$ whose silhouette has size $s$) has no $ab$-cycle of length in $[2,n^\alpha]$. Then
$$\probap'_\alpha(n) = \sum_{s = n-3n^{2/3}}^n \P_n\left(\cyc_n(s)\right)\ \probap'_\alpha(n,s),$$
where $\P_{n}$ is the uniform probability on $\cyc_n$.

Finally, for $\Gamma$ taken uniformly at random in $\cyc_n(s)$ let $\probaq_\alpha(n,s)$ be the probability that both $\core(\Gamma)$ and $\Gamma$ have no $ab$-cycle of size in $[2,n^\alpha]$, and $\probaq'_\alpha(n,s)$ the probability that $\core(\Gamma)$ has such a cycle but $\Gamma$ does not. Then
$$\probap'_\alpha(n,s) \le \probaq_\alpha(n,s) + \probaq'_\alpha(n,s).$$

Observe that $\probaq_\alpha(n,s)$ is at most equal to the probability that $\Gamma$ or $\core(\Gamma)$ has no $ab$-cycle of size in $[2,n^\alpha]$. 
Theorem~\ref{thm:uniformity} shows that the latter is equal to the probability that a size $s$ silhouette graph has no $ab$-cycle of size in $[2,s^\alpha]$, which is $\O(s^{-\alpha})$ according to Proposition~\ref{pro:small cycles}. Thus, there exist positive constants $C$ (independent of $s$ and $n$) such that
$$\probaq_\alpha(n,s) \le Cs^{-\alpha} \le C(n-3n^{\frac23})^{-\alpha}$$

Turning our attention to $\probaq'_\alpha(n,s)$, we note that
$$\probaq'_\alpha(n,s) = \sum \P_{n,s}\left(\core(\Gamma) = \Delta\right)\ \probaq_\Delta,$$
where $\P_{n,s}$ is the uniform probability on $\cyc_n(s)$, $\probaq_\Delta$ is the probability that a graph in $\cyc_n(s)$ having silhouette $\Delta$ has no $ab$-cycle of size in $[2,n^\alpha]$, and the sum is taken over all size $s$ silhouette graphs $\Delta$ who do have an $ab$-cycle of such a size.

Let $\Delta = \core(\Gamma)$. By Proposition~\ref{prop: uniqueness silhouette}, $\Delta$ can be obtained from $\Gamma$ after a maximal sequence of of $\lambda_3$-, $\lambda_{2,1}$-, $\lambda_{2,2}$- and $\kappa_{3}$-moves (since $s \ge n-3n^{\frac23} > 2$). In particular, the $a$-edges of $\Delta$ that were not already edges in $\Gamma$ arose during a $\kappa_3$-move, which deleted two vertices, two $a$-edges and a $b$-edge, and added a new $a$-edge.

Since $\Gamma$ has at most $3n^{\frac{2}{3}}$ vertices more that than $\Delta$, the path from $\Gamma$ to $\Delta$ has at most $\frac32n^{\frac{2}{3}}$ $\kappa_3$-moves, and the majority of $a$-edges of $\Delta$ are also $a$-edges of $\Gamma$. More precisely, at most $\frac32n^{\frac{2}{3}}$ of the $\frac s2$ $a$-edges of $\Delta$ are not $a$-edges of $\Gamma$.
Thus the probability that an $a$-edge of $\Delta$ fails to be an edge of $\Gamma$ is at most
$$ \frac{3n^{\frac23}}{s}\enspace\le\enspace \frac{3n^{\frac23}}{n-3n^{\frac23}} \enspace=\enspace 3n^{-\frac13}\left(1+\O\left(n^{-\frac13}\right)\right).$$

For each of the silhouette graphs $\Delta$ under consideration (of size $s$, with an $ab$-cycle of size in $[2,n^\alpha]$), fix an $ab$-cycle in $\Delta$ of size, say, $\lambda$. The probability that at least one of the $\lambda$ $a$-edges in this $ab$-cycle is not an edge in $\Gamma$ (that is: that this cycle does not exist in $\Gamma$), is bounded above by
$$3\lambda n^{-\frac13}\left(1+\O\left(n^{-\frac13}\right)\right) \enspace\le\enspace 3n^{\alpha-\frac13}\left(1+\O\left(n^{-\frac13}\right)\right).$$
As a result,
$\probaq_\Delta \enspace\le\enspace 3n^{\alpha-\frac13}\left(1+\O\left(n^{-\frac13}\right)\right), 
\probaq'_\alpha(n,s) \enspace\le\enspace 3n^{\alpha-1/3}\left(1+\O\left(n^{-1/3}\right)\right)$, and 
$$\probap'_\alpha(n,s) \enspace\le\enspace \probaq_\alpha(n,s) + \probaq'_\alpha(n,s) \enspace\le\enspace C'n^{-\alpha} + 3n^{\alpha-1/3}\left(1+\O\left(n^{-1/3}\right)\right).$$
Since $\alpha < \frac16$, $\alpha - \frac13 < -\alpha$ and there exists a constant $C''  > 0$ such that $\probap'_\alpha(n,s) \le C''n^{-\alpha}$, and hence $\probap'_\alpha(n) \le C''n^{-\alpha}$ which, in turn, yields $\probap_\alpha(n) = \O(n^{-\alpha})$.

The same proof holds for $\PSL$-reduced graphs, reasoning within $\rooted_n$ instead of $\cyc_n$, and using Proposition~\ref{prop: large silhouette rooted} and Theorem~\ref{thm:uniformity rooted} instead of Proposition~\ref{prop: large silhouette} and Theorem~\ref{thm:uniformity}.
\eopo

\subsection{Proof of \protect{Proposition~\ref{pro:small cycles}}}\label{sec: cycles in silhouette}

In this section, we give a proof of Proposition \ref{pro:small cycles} which states that, with high probability, in a finite group of permutations generated by a pair of permutations without fixed points $(\sigma_2,\sigma_3)$, of order 2 and 3 respectively, the composition $\sigma_2\sigma_3$ admits orbits of a certain, relatively small size. This type of result, on the composition of two randomly chosen mappings, is notoriously difficult to obtain and is therefore interesting {\it per se}.

Let $0 < \alpha < \frac16$ and let $n > 0$ be a multiple of 6. We want to show that a random size $n$ silhouette graph fails to have an $ab$-cycle of size at least $2$ and at most $n^\alpha$ with probability $\O(n^{-\alpha})$.

Let $\calG$ be the set of labeled graphs that are disjoint unions of silhouette graphs. Recall that the set of $n$-vertex elements of $\calG$ ($n$ a multiple of 6) is in bijection with the set of pairs $(\sigma_2,\sigma_3)$ of fixpoint-free permutations of $[n]$, the first of order 2 (the $a$-edges) and the second of order 3 (the $b$-edges).

If we fix $\sigma_3$, the set $\calG_{\sigma_3}$ of elements of $\calG$ characterized by pairs of the form $(\sigma_2,\sigma_3)$ has cardinality the number of fixpoint-free, order 2 permutations on $[n]$, namely $(n-1)!!$, where the double factorial of an odd integer $q$ is given by $q!! = q(q-2)(q-4)\dots 1$. Since this value does not depend on $\sigma_3$, it is sufficient to establish the restriction of Proposition~\ref{pro:small cycles} to an arbitrary $\calG_{\sigma_3}$, provided that the constants in the $\O$-notation do not depend on $\sigma_3$.

For convenience, we fix $\sigma_3 = (1\,2\,3)(4\,5\,6)\dots(n-2\,n-1\,n)$, and we write $\calG_n$ instead of $\calG_{\sigma_3}$.
We now concentrate on a particular set of $ab$-cycles. Say that an $ab$-cycle in a graph $G \in \calG$ is \emph{simple} if it visits at most one vertex in each $b$-triangle in $G$. Equivalently, walking along the corresponding $(ab)^m$-labeled path does not require traveling through an $a$-edge in both directions. 

Now let $M = \lfloor n^\alpha\rfloor$ and $\MM = \{2,\dots,M\}$. In the context of this proof, we say that a set is \emph{small} if its cardinality is in $\MM$. Then a \emph{small} $ab$-cycle in an element of $\calG_n$ visits a \emph{small} set of $b$-triangles. We need to prove that, only with probability $\O(n^{-\alpha})$, a random size $n$ silhouette graph fails to have a small simple $ab$-cycle.

Let $\calC$ be the set of graphs in $\calG_n$ with at least a small simple $ab$-cycle. Let also $\calS$ be the set of all small sets of $b$-triangles in a graph of $\calG_n$ (recall that all the graphs of $\calG_n$ have the same $b$-triangles). If $I \in \calS$, let $\calC_I$ be the set of graphs in $\calG_n$ containing a simple $ab$-cycle visiting exactly the $b$-triangles in $I$. Then
$
\calC = \bigcup_{I\in\calS} \calC_I.
$

Finally, let $\ocalC$ denote the complement of $\calC$ in $\calG_n$ (the elements of $\calG_n$ without a small simple $ab$-cycle). We want to show that $\frac{|\ocalC|}{|\calG_n|} = \frac{|\ocalC|}{(n-1)!!}$ is of the form $\O(n^{-\alpha})$.
By the inclusion-exclusion principle, we have
$$|\ocalC| = |\calG_n| - \sum_{I\in\calS}|\calC_I| + \sum_{\substack{\{I_1,I_2\}\subseteq\calS\\ I_1\ne I_2}} |\calC_{I_1}\cap \calC_{I_2}| - \ldots
= \sum_{\calI\subseteq\calS}(-1)^{|\calI|}\left| \bigcap_{I\in\calI}\calC_I\right|
= \sum_k (-1)^k \sum_{\substack{\calI\subseteq S\\|\calI| = k}}\left| \bigcap_{I\in\calI}\calC_I\right|.$$
Truncating the inclusion-exclusion formula on even or odd cardinalities for $\calI$, yields upper and lower bounds for $|\ocalC|$. For any $\kappa \geq 0$ we have
\begin{equation}\label{eq: inclusion-exclusion}
\sum_{k=0}^{2\kappa+1} (-1)^k \sum_{\substack{\calI\subseteq S\\|\calI| = k}}\left| \bigcap_{I\in\calI}\calC_I\right|
\enspace\leq\enspace |\ocalC| \enspace\leq\enspace
 \sum_{k=0}^{2\kappa} (-1)^k \sum_{\substack{\calI\subseteq S\\|\calI| = k}}\left| \bigcap_{I\in\calI}\calC_I\right|
\end{equation}

It turns out to be more convenient to work with tuples of $b$-triangles rather than sets. If $\calI=\{I_1,\ldots,I_k\}$ is a small set of $b$-triangles with cardinality $k$, there are $k!$ tuples
$\bm J = (J_1,\ldots,J_k)\in\calS^k$ such that $\{J_1,\ldots,J_k\}=\calI$, so
\begin{equation}
\left| \bigcap_{I\in\calI}\calC_I\right| = \frac1{k!}\sum_{\substack{(J_1,\ldots,J_k)\in\calS^k\\\{J_1,\ldots,J_k\}=\calI}} \left| \calC_{J_1}\cap\ldots\cap\calC_{J_k}\right|.
\label{eq: given set of triangles}
\end{equation}

Say that two simple $ab$-cycles \emph{overlap} if they visit a same $b$-triangle (of course, on different vertices). We now distinguish the tuples $\bm J = (J_1,\ldots,J_k)$ in Equation~\eqref{eq: given set of triangles} according to the cardinality of their components and to their overlaps. More precisely, if $\bm d = (d_1,\ldots,d_k) \in \MM^k$, we let
\begin{align*}
\nooverlap(\bm d) &= \left\{ \bm J = (J_1,\ldots,J_k)\in\calS^k \mid \forall i,\ |J_i|=d_i\text{ and the $J_i$ are pairwise disjoint}\right\} \\
\overlap(\bm d) &= \Big\{ \bm J = (J_1,\ldots,J_k)\in\calS^k \mid \forall i,\ |J_i|=d_i\text{ and  the $J_i$ are pairwise distinct} \\
&\hskip 8cm\text{but not pairwise disjoint}\Big\}.
\end{align*}

Returning to the summands in the estimation of the inclusion-exclusion bounds (Equation~\eqref{eq: inclusion-exclusion}), we now have, for every $k$,
\begin{equation}\label{eq: Ak Bk}
\sum_{\substack{\calI\subseteq S\\|\calI| = k}}\left| \bigcap_{I\in\calI}\calC_I\right|
\enspace=\enspace \frac1{k!}\left(\underbrace{\sum_{\bm d\in\MM^k}\sum_{\bm J\in\nooverlap(\bm d)}\left| \bigcap_{i=1}^k\calC_{J_i}\right|}_{A_k}+\underbrace{\sum_{\bm d\in\MM^k}\sum_{\bm J\in\overlap(\bm d)}\left| \bigcap_{i=1}^k\calC_{J_i}\right|}_{B_k}\right)
\end{equation}
We now study successively the quantities $A_k$ and $B_k$ in Equation~\eqref{eq: Ak Bk}. If $q\ge 1$, we let $H_q$ denote the partial sum $\sum_{i=1}^q\frac1i$ of the harmonic series. It is well known that $H_q = \log q + \gamma + o(1)$, where $\gamma$ is Euler's constant.

\begin{lemma}\label{lemma:Ak}
Let $n$ be a positive multiple of $6$. Let $0 < \alpha < \frac12$, $M = \lfloor n^\alpha\rfloor$, $0 < \beta < \frac12-\alpha$ and $1 \le k \le n^\beta$. Finally, let $\delta=\alpha+\beta$ and let $A_k$ be as in Equation~\eqref{eq: Ak Bk}. Then 
\begin{equation*}
\frac1{(n-1)!!}A_k = (H_M-1)^k\left(1 + \O\left(n^{2\delta-1}\right)\right),
\end{equation*}
uniformly in $k$ (that is: the constants intervening in the $\O$ notation do not depend on $n$ or $k$).
\end{lemma}

\proof
Let $\bm d = ( d_1,\dots, d_k) \in\MM^k$, $ d= d_1+\ldots+ d_k$ and $\bm J = (J_1,\ldots,J_k)\in\nooverlap(\bm d)$. To construct a graph in $\bigcap_i\calC_{J_i}$, that is, a graph in $\calG_n$ with (non-overlapping) simple $ab$-cycles over the sets of $b$-triangles $J_1, \dots, J_k$ respectively, we must 
\begin{itemize}
\item select for each $b$-triangle in $J_1, \dots, J_k$ a vertex belonging to the collection of simple $ab$-cycles, {\it i.e.} a vertex reached after reading an occurrence of $ab$ in an $ab$-cycle; there are 3 possibilities for each $b$-triangle, and therefore a total of $3^d$ choices;
\item cyclically order the triangles in each $J_i$; there are $( d_1-1)!\dots( d_k-1)!$ possibilities to do so; note that this second step fully determines the $a$-edges adjacent to the $ d$ vertices chosen in the first step;
\item choose the missing $a$-edges arbitrarily: they connect the $n-2 d$ vertices not yet adjacent to an $a$-edge, and there are $(n-2 d-1)!!$ ways to do so.
\end{itemize}

Thus, for every $\bm J = (J_1,\ldots,J_k)\in\nooverlap(\bm d)$, there are $3^ d( d_1-1)!\cdots( d_k-1)!(n-2 d-1)!!$ graphs in $\bigcap_{i=1}^k\calC_{J_i}$ (which is independent of the choice of $\bm J$ in $\nooverlap(\bm d)$).

Moreover, since a graph in $\calG_n$ has $\frac n3$ $b$-triangles and the components of a tuple in $\nooverlap(\bm d)$ are pairwise disjoint, we have $|\nooverlap(\bm d)| = \binom{n/3}{ d_1,\ldots, d_k,n/3- d}$. Thus
$$\sum_{\bm J\in\nooverlap(\bm d)}\left| \bigcap_{i=1}^k\calC_{J_i}\right|
= \binom{n/3}{ d_1,\ldots, d_k,n/3- d}3^d( d_1-1)!\cdots( d_k-1)!(n-2 d-1)!!.$$
Note that $\frac{(n-2  d-1)!!}{(n-1)!!}= \prod_{i=0}^{ d-1}\frac{1}{n-1-2i}$ and $\frac{(n/3)! 3^{ d}}{(n/3- d)!}=\prod_{i=0}^{ d-1} (n-3i)$. Therefore
$$\frac1{(n-1)!!}\ \sum_{\bm J\in\nooverlap(\bm d)}\left| \bigcap_{i=1}^k\calC_{J_i}\right| \enspace=\enspace  \frac{Q_d}{ d_1\cdots d_k},\quad\text{ with } Q_d = \prod_{i=0}^{d-1}\frac{n-3i}{n-1-2i}.$$
Now observe that, for $n$ sufficiently large,
\begin{align}
Q_d &\enspace\leq\enspace \frac{n}{n-1}\quad\text{and} \label{eq: upper bound of Qell}\\
Q_d &\enspace=\enspace  \prod_{i=0}^{d-1}\left(1-\frac{i-1}{n-1-2i} \right) \enspace\ge\enspace \left(1-\frac{d}{n-2d}\right)^d \enspace=\enspace \exp\left(d\log\left(1-\frac{d}{n-2d}\right)\right) \nonumber \\
&\enspace\ge\enspace 1 + d\log\left(1-\frac{d}{n-2d}\right).\label{eq: lower bound of Qell}
\end{align}
The lower bound in Equation \eqref{eq: lower bound of Qell} for $Q_d$ is a decreasing function of $d$ and the possible values of $d$ satisfy  $1\le d \le kM \le n^{\alpha+\beta} = n^\delta$, so
\begin{equation}\label{eq: U and L}
1 + n^\delta\log\left(1-\frac{n^\delta}{n-2n^\delta}\right) \enspace\le\enspace Q_d \enspace\le\enspace \frac n{n-1}.
\end{equation}
Let $U$ and $L$ be the upper and lower bounds of $Q_d$ in Equation~\eqref{eq: U and L}. Note that $\sum_{d\in\MM^k}\frac1{d_1\cdots d_k} = (H_M-1)^k$. Therefore
$$L (H_M-1)^k \enspace=\enspace L \sum_{\bm d\in\MM^k}\frac1{d_1\dots d_k} \enspace\le\enspace \frac1{(n-1)!!}A_k \enspace\le\enspace  U \sum_{\bm d\in\MM^k}\frac1{d_1\dots d_k} \enspace=\enspace U(H_M-1)^k,$$
which concludes the proof of Lemma~\ref{lemma:Ak} since both $U$ and $L$ equal $1 + \O\left(n^{2\delta-1}\right)$.
\eop

\begin{lemma}\label{lemma:Bk}
 Let $n$ be a positive multiple of $6$. Let $0 < \alpha < \frac16$, $M = \lfloor n^\alpha\rfloor$, $0 < \beta < \frac16-\alpha$ and $1 \le k \le n^\beta$. Finally, let $\delta = \alpha+\beta$ and let $B_k$ be as in Equation~\eqref{eq: Ak Bk}. Then
\begin{equation*}
\frac1{(n-1)!!} B_k \leq \O\left(n^{6\delta-1}\right) (H_M-1)^k
\end{equation*}
uniformly in $k$.
\end{lemma}

\proof
Let $\bm d = (d_1,\dots,d_k) \in\MM^k$ and $d=d_1+\ldots+d_k$. Since every vertex of a silhouette graph occurs in exactly one $ab$-cycle, a $b$-triangle can occur in at most 3 $ab$-cycles. Let $\bm t =(t_1,t_2,t_3)$ be a tuple of non-negative integers such that $d = t_1+2t_2+3t_3$. We denote by $\overlap_k(\bm d;\bm t)$ the set of elements $\bm J = (J_1,\ldots,J_k)\in\overlap_k(\bm d)$ such that there are $t_1$ (respectively, $t_2$, $t_3$) $b$-triangles belonging to exactly 1 (respectively, 2, 3) of the components of $\bm J$. We talk of $b$-triangles \emph{of type} 1 (respectively, 2, 3).

\begin{claim}\label{claim: even}
If $\bm J \in \overlap_k(\bm d;\bm t)$ and $\bigcap_{i=1}^k\calC_{J_i} \ne \emptyset$, then $t_2+t_3$ is even.
\end{claim}

\proofof{Claim~\ref{claim: even}}
Let $G \in \bigcap_{i=1}^k\calC_{J_i}$. Consider the $a$-edges occurring in the $k$ small simple $ab$-cycles visiting, respectively, $J_1,\dots, J_k$. It is convenient at this point to think of the $a$-edges in our $ab$-cycles as matched pairs of half-edges, which we denote $(T,e)$, where $e$ is an $a$-edge and $T$ is a $b$-triangle adjacent to $e$. Since no $a$-edge in a simple $ab$-cycle may connect vertices from the same $b$-triangle, each such $a$-edge $e$ corresponds to a pair of distinct, matched half-edges $(T,e)$ and $(T',e)$.

Even though $a$-edges are undirected (or can be traversed in both directions), considering one in a simple $ab$-cycle uniquely defines a direction and we can talk of matched outgoing and incoming half-edges along a simple $ab$-cycle. Clearly, an $a$-edge is used in only one direction in a given $ab$-cycle, but it can be used in different directions by distinct $ab$-cycles. So we say that a half-edge $(T,e)$ such that $e$ occurs in the union of the $k$ $ab$-cycles under consideration is \emph{outgoing} (respectively, \emph{incoming}, \emph{$2$-way}) if it only occurs as outgoing (respectively, it only occurs as incoming, it occurs both as outgoing and incoming).

If $T$ is a $b$-triangle in $\bigcup J_i$, occurring in just one of the $k$ small simple $ab$-cycles, then $T$ is a component of exactly 1 incoming and 1 outgoing half-edges; if $T$ occurs in two of these $ab$-cycles, it is a component of 1 incoming, 1 outgoing and 1 $2$-way half-edges; finally, if $T$ occurs in three $ab$-cycles, it is a component of 3 $2$-way half-edges.

The result follows since a $2$-way half-edge must be matched with another, distinct, $2$-way half-edge, and we have $t_2+3t_3$ such half-edges.
\eopo

Let $\bm J = (J_1,\ldots,J_k) \in \overlap_k(\bm d;\bm t)$. Then $|\bigcup_i J_i| = t_1+t_2+t_3$ whereas $d = \sum_i|J_i| = t_1+2t_2+3t_3$.  The overlaps between the $J_i$ determine which $b$-triangles in $\bigcup_i J_i$ are of type 1, 2 or 3. To construct an $n$-vertex graph in $\bigcap_i\calC_{J_i}$, that is, a graph in $\calG_n$ with (overlapping) simple $ab$-cycles over the sets of $b$-triangles $J_1, \dots, J_k$ respectively, we must
\begin{itemize}
\item select for each $b$-triangle of type 1 a vertex belonging to the collection of simple $ab$-cycles (recall that it is a vertex reached after reading an occurrence of $ab$ in an $ab$-cycle); there are $3^{t_1}$ choices. Similarly, for each $b$-triangle of type 2, select two vertices belonging to the collection of simple $ab$-cycles; there are $3^{t_2}$ choices. Note that every vertex of a type 3 $b$-triangle belongs to the collection of simple $ab$-cycles. As discussed in the proof of Claim~\ref{claim: even}, the choice of these $t_1+2t_2+3t_3$ vertices in the $ab$-cycles implies the presence of $2t_1+3t_2+3t_3$ vertices along the corresponding $(ab)^{d_i}$-labeled cycles.
\item cyclically order the triangles in each $J_i$; there are $(d_1-1)!\dots(d_k-1)!$ possibilities to do so; note that this second step fully determines the $a$-edges adjacent to the $2t_1 + 3t_2 + 3t_3$ vertices determined in the first step;
\item choose the missing $a$-edges arbitrarily: they connect the $n-(2t_1 + 3t_2 + 3t_3)$ vertices not yet selected, and there are $(n-(2t_1 + 3t_2 + 3t_3)-1)!!$ ways to do so.
\end{itemize}

Thus, for every $\bm J = (J_1,\ldots,J_k) \in \overlap_k(\bm d;\bm t)$, we have
$$\left|\bigcap_{i=1}^k\calC_{J_i}\right| \enspace=\enspace 3^{t_1+t_2}(d_1-1)!\cdots(d_k-1)!(n-2t_1 - 3t_2 - 3t_3 -1)!!.$$
Now, to construct a tuple in $\overlap_k(\bm d;\bm t)$, we must 
\begin{itemize}
\item choose $t_i$ $b$-triangles of type $i$ ($i = 1, 2, 3$); there are
$\binom{n/3}{t_1,t_2,t_3,n/3-t_1-t_2-t_3}$ choices;
\item allocate the $t_1+t_2+t_3$ selected $b$-triangles to sets $J_1,\dots,J_k$ respecting multiplicities and the required cardinality of the $J_i$; there are at most  $\binom{t_1+2t_2+3t_3}{d_1,\ldots,d_k} = \binom{d}{d_1,\ldots,d_k}$ choices (this is an upper bound as some choices may be unrealizable or produce cycles that are not simple).
\end{itemize}

Therefore we have
$$\left|\overlap_k(\bm d;\bm t)\right| \enspace\le\enspace \binom{n/3}{t_1,t_2,t_3,n/3-t_1-t_2-t_3}\ \binom{d}{d_1,\ldots,d_k}$$
and $\sum_{\bm J\in\overlap(\bm d;\bm t)}\left| \bigcap_{i=1}^k\calC_{J_i}\right|$ is bounded above by
$$\binom{n/3}{t_1,t_2,t_3,n/3-t_1-t_2-t_3}\ \binom{d}{d_1,\ldots,d_k}3^{t_1+t_2}(d_1-1)!\cdots(d_k-1)!(n-2t_1 - 3t_2 - 3t_3 -1)!!.$$
Dividing by $(n-1)!!$, we get
$$\frac{\sum_{\bm J\in\overlap(\bm d;\bm t)}\left| \bigcap_{i=1}^k\calC_{J_i}\right|}{(n-1)!!} \enspace\le\enspace \frac{d !}{d_1\cdots d_k\ t_1!t_2!t_3!\  3^{t_3}}\ 
\frac{n(n-3)\cdots(n-3(t_1+t_2+t_3-1))}{(n-1)(n-3)\cdots(n-2t_1-3t_2-3t_3+1)}.$$
Let $t=t_1+t_2+t_3$, so that $t_1=d - 2t_2-3t_3$. Recall that $Q_t =  \prod_{i=0}^{t-1}\frac{n-3i}{n-1-2i}\leq \frac{n}{n-1}$ (Equation~\eqref{eq: upper bound of Qell}). Then, as $t_2!t_3!3^{t_3}\geq1$,
\begin{align}
\frac{\sum_{\bm J\in\overlap(\bm d;\bm t)}\left| \bigcap_{i=1}^k\calC_{J_i}\right|}{(n-1)!!}
& \enspace\leq\enspace \frac1{d_1\cdots d_k}\cdot \frac{d!}{(d-2t_2-3t_3)!}
\cdot\frac{Q_t}{\prod_{i=t}^{t_1+\frac32(t_2+t_3)-1}(n-1-2i)} \nonumber\\
& \enspace\leq\enspace \frac{n\ d^{2t_2+3t_3}}{(n-1)\ d_1\cdots d_k\ (n-2t_1-3t_2-3t_3+1)^{\frac12(t_2+t_3)}} \label{eq: step in Bk}
\end{align}
Since $d = t_1 + 2t_2 + 3t_3$ and $d \leq kM \leq n^{\alpha+\beta} = n^\delta$, we have $2t_1+3t_2+3t_3 \leq 2d \leq 2n^\delta$, yielding
\[
\frac{d^{2t_2+3t_3}}{(n-2t_1-3t_2-3t_3+1)^{\frac12(t_2+t_3)}} \leq 
\frac{n^{\delta(2t_2+3t_3)}}{(n-2n^\delta)^{\frac12(t_2+t_3)}}.
\]
Moreover, $\frac12(t_2+t_3)\leq \frac14 d \leq \frac14 n^\delta$ and we have
\begin{align*}
(n-2n^\delta)^{\frac12(t_2+t_3)}
& = n^{\frac12(t_2+t_3)}\exp\left(\left(\frac12(t_2+t_3)\right)\log\left(1-2n^{\delta-1}\right)\right)\\
& \geq n^{\frac12(t_2+t_3)}\exp\left(\frac14 n^\delta \log\left(1-2n^{\delta-1}\right)\right)
\end{align*}
The term under the exponential is
$\frac14 n^\delta \log\left(1-2n^{\delta-1}\right) \sim -\frac12 n^{3\delta-1},$
which tends to 0 (since $\delta < \frac16 < \frac13$). Therefore, for $n$ sufficiently large, 
\[
(n-2n^\delta)^{\frac12(t_2+t_3)} \geq \frac12 n^{\frac12(t_2+t_3)}.
\]
Going back to Equation~\eqref{eq: step in Bk}, we have
\begin{equation}\label{eq: 2nd step in Bk}
\frac{\sum_{\bm J\in\overlap(\bm d;\bm t)}\left| \bigcap_{i=1}^k\calC_{J_i}\right|}{(n-1)!!} \leq  \frac{1}{d_1\cdots d_k} \frac{n^{\delta(2t_2+3t_3)}}{n^{\frac12(t_2+t_3)}} = \frac{1}{d_1\cdots d_k} \ n^{\delta(2t_2+3t_3)-\frac12(t_2+t_3)}\\
\end{equation}
Equation~\eqref{eq: 2nd step in Bk} must be summed over $\bm t$ and $\bm d$. Let us first fix $\bm d$ and recall that $\bm t=(t_1,t_2,t_3)$ with $(t_2,t_3)\neq (0,0)$ and $t_2+t_3$ even (Claim~\ref{claim: even}).

Consider first the $\bm t$ of the form $(t_1,t_2,0) = (t_1,2s,0)$ ($s\ge 1$). The corresponding subsum of powers of $n$ is bounded above (for $n\ge 2$) by
$$\sum_{s\ge 1}\left(n^{2\delta-\frac12}\right)^{2s} = \frac{n^{4\delta -1}}{1 - n^{4\delta -1}} \le 2 n^{4\delta-1}.$$
Similarly, the subsum of powers of $n$ corresponding to the $\bm t$ of the form $(t_1,0,t_3)$ is bounded above by $2 n^{6\delta-1}$.

Next, the subsum corresponding to the $\bm t$ where both $t_2$ and $t_3$ are non-zero is bounded above by
$$\sum_{t_2\ge 1}\sum_{t_3\ge 1} \left(n^{2\delta-\frac12}\right)^{t_2}\left(n^{3\delta-\frac12}\right)^{t_3} \le \left(\sum_{t_2\ge 1}\left(n^{2\delta-\frac12}\right)^{t_2}\right)\ \left(\sum_{t_3\ge 1}\left(n^{3\delta-\frac12}\right)^{t_3}\right) \le 2 n^{5\delta-1}$$
It follows that, for a fixed tuple $\bm d$,
\[
\frac{\sum_{\bm J\in\overlap(\bm d;\bm t)}\left| \bigcap_{i=1}^k\calC_{J_i}\right|}{(n-1)!!} 
\leq \frac{2(n^{4\delta-1} + n^{5\delta-1} + n^{6\delta-1})}{d_1\cdots d_k}.
\]
Now, summing over $\bm d$, we get
\begin{align*}
\frac{\sum_{\bm J\in\overlap(\bm d;\bm t)}\left| \bigcap_{i=1}^k\calC_{J_i}\right|}{(n-1)!!} 
&\leq 2(n^{4\delta-1} + n^{5\delta-1} + n^{6\delta-1})\sum_{\bm d}\frac1{d_1\cdots d_k} \\
&\le 2(n^{4\delta-1} + n^{5\delta-1} + n^{6\delta-1})(H_M-1)^k
\end{align*}
Since $n^{4\delta-1} + n^{5\delta-1} + n^{6\delta-1} = \O\left(n^{6\delta-1}\right)$, we have, for $n$ large enough,
\[
\frac{\sum_{\bm J\in\overlap(\bm d;\bm t)}\left| \bigcap_{i=1}^k\calC_{J_i}\right|}{(n-1)!!} = (H_M-1)^k \O(n^{6\delta-1}),
\]
thus concluding the proof of Lemma~\ref{lemma:Bk}.
\eop

We can now conclude the proof of Proposition~\ref{pro:small cycles}. 
By Equations~\eqref{eq: inclusion-exclusion} and~\eqref{eq: Ak Bk}, we want to show that, for $\kappa=\lfloor n^\beta\rfloor$ and $\kappa=\lfloor n^\beta\rfloor -1$,
\[
\frac1{(n-1)!!}\ \sum_{k=0}^{\kappa}(-1)^k\sum_{\substack{\calI\subseteq S\\|\calI| = k}}\left| \bigcap_{I\in\calI}\calC_I\right|
= \sum_{k=0}^{\kappa}\frac{(-1)^k}{k!}\frac{A_k}{(n-1)!!} + \sum_{k=0}^{\kappa}\frac{(-1)^k}{k!}\frac{B_k}{(n-1)!!}
\]
is $\O(n^{-\alpha})$. 
By Lemma~\ref{lemma:Ak}, the absolute value of the first sum is bounded above by
$$(1+\O(n^{2\delta-1})) \exp(-(H_M-1)) = e^{1-\gamma}n^{-\alpha}(1+\O(n^{2\delta-1})).$$
And by Lemma~\ref{lemma:Bk}, the absolute value of the second sum is bounded above by
$$\exp(-(H_M-1))\ \O(n^{6\delta-1}) = e^{1-\gamma}n^{-\alpha} \O(n^{6\delta-1}).$$
Thus the whole sum is $\O(n^{-\alpha})$, establishing the expected bound for disjoint unions of silhouette graphs.

Such a union, of size $n$ (a multiple of 6), is connected (and hence silhouette) with probability $1 - \frac56n\inv + o(n\inv)$ by \cite[Proof of Proposition 8.18]{2021:BassinoNicaudWeil}. Thus the probability that a silhouette graph has no small simple $ab$-cycle is, again, $\O(n^{-\alpha})$.

{\small\bibliographystyle{abbrv}

\begin{thebibliography}{10}

\bibitem{1998:Arzhantseva}
G.~N. Arzhantseva.
\newblock Generic properties of finitely presented groups and {H}owson's
  theorem.
\newblock {\em Comm. Algebra}, 26(11):3783--3792, 1998.

\bibitem{1996:ArzhantsevaOlshanskii}
G.~N. Arzhantseva and A.~Y. Ol'shanski{\u \i}.
\newblock Generality of the class of groups in which subgroups with a lesser
  number of generators are free.
\newblock {\em Mat. Zametki}, 59(4):489--496, 1996.

\bibitem{2013:BassinoMartinoNicaud}
F.~Bassino, A.~Martino, C.~Nicaud, E.~Ventura, and P.~Weil.
\newblock Statistical properties of subgroups of free groups.
\newblock {\em Random Struct. Algorithms}, 42(3):349--373, 2013.

\bibitem{2008:BassinoNicaudWeil}
F.~Bassino, C.~Nicaud, and P.~Weil.
\newblock Random generation of finitely generated subgroups of a free group.
\newblock {\em Internat. J. Algebra Comput.}, 18(2):375--405, 2008.

\bibitem{2016:BassinoNicaudWeilCM}
F.~Bassino, C.~Nicaud, and P.~Weil.
\newblock Generic properties of subgroups of free groups and finite
  presentations.
\newblock In {\em Algebra and computer science}, volume 677 of {\em Contemp.
  Math.}, pages 1--43. Amer. Math. Soc., Providence, RI, 2016.

\bibitem{2016:BassinoNicaudWeil}
F.~Bassino, C.~Nicaud, and P.~Weil.
\newblock On the genericity of {W}hitehead minimality.
\newblock {\em J. Group Theory}, 19(1):137--159, 2016.

\bibitem{2021:BassinoNicaudWeil}
F.~Bassino, C.~Nicaud, and P.~Weil.
\newblock Statistics of subgroups of the modular group.
\newblock {\em {Int. J. Algebra Comput.}}, 31(8):1691--1751, 2021.

\bibitem{2023:BassinoNicaudWeil}
F.~Bassino, C.~Nicaud, and P.~Weil.
\newblock Random generation of subgroups of the modular group with a fixed
  isomorphism type.
\newblock \texttt{arXiv:2310.18923}, 2023.

\bibitem{2021:BudzinskiCurienPetri}
T.~Budzinski, N.~Curien, and B.~Petri.
\newblock The diameter of random {B}elyi surfaces.
\newblock {\em Algebraic \& Geometric Topology}, 21:2929--2957, 2021.

\bibitem{2009:FlajoletSedgewick}
P.~Flajolet and R.~Sedgewick.
\newblock {\em Analytic combinatorics}.
\newblock Cambridge University Press, 2009.

\bibitem{2006:Gamburd}
A.~Gamburd.
\newblock Poisson-{D}irichlet distribution for random {B}elyi surfaces.
\newblock {\em Ann. Probab.}, 34(5):1827--1848, 2006.

\bibitem{1991:GerstenShort}
S.~M. Gersten and H.~B. Short.
\newblock Rational subgroups of biautomatic groups.
\newblock {\em Ann. of Math.}, 134(1):125--158, 1991.

\bibitem{2010:GilmanMiasnikovOsin}
R.~Gilman, A.~Miasnikov, and D.~Osin.
\newblock Exponentially generic subsets of groups.
\newblock {\em Illinois J. Math.}, 54(1):371--388, 2010.

\bibitem{1996:Gitik}
R.~Gitik.
\newblock Nielsen generating sets and quasiconvexity of subgroups.
\newblock {\em J. Pure Appl. Algebra}, 112(3):287--292, 1996.

\bibitem{1990:JacksonVisentin}
D.~M. Jackson and T.~I. Visentin.
\newblock A character-theoretic approach to embeddings of rooted maps in an
  orientable surface of given genus.
\newblock {\em Trans. Am. Math. Soc.}, 322(1):343--363, 1990.

\bibitem{2002:Jitsukawa}
T.~Jitsukawa.
\newblock Malnormal subgroups of free groups.
\newblock In {\em Computational and statistical group theory ({L}as {V}egas,
  {NV}/{H}oboken, {NJ}, 2001)}, volume 298 of {\em Contemp. Math.}, pages
  83--95. Amer. Math. Soc., Providence, RI, 2002.

\bibitem{1996:Kapovich}
I.~Kapovich.
\newblock Detecting quasiconvexity: algorithmic aspects.
\newblock In {\em Geometric and computational perspectives on infinite groups
  ({M}inneapolis, {MN} and {N}ew {B}runswick, {NJ}, 1994)}, volume~25 of {\em
  DIMACS Ser. Discrete Math. Theoret. Comput. Sci.}, pages 91--99. Amer. Math.
  Soc., Providence, RI, 1996.

\bibitem{2017:KharlampovichMiasnikovWeil}
O.~Kharlampovich, A.~Miasnikov, and P.~Weil.
\newblock Stallings graphs for quasi-convex subgroups.
\newblock {\em {J. Algebra}}, 488:442--483, 2017.

\bibitem{1977:LyndonSchupp}
R.~C. Lyndon and P.~E. Schupp.
\newblock {\em Combinatorial group theory}.
\newblock Springer-Verlag, 1977.

\bibitem{1976:MagnusKarrassSolitar}
W.~Magnus, A.~Karrass, and D.~Solitar.
\newblock {\em Combinatorial group theory}.
\newblock Dover, 1976.

\bibitem{2019:MaherSisto}
J.~Maher and A.~Sisto.
\newblock Random subgroups of acylindrically hyperbolic groups and hyperbolic
  embeddings.
\newblock {\em Int. Math. Res. Not. IMRN}, 13:3941--3980, 2019.

\bibitem{2007:Markus-Epstein}
L.~Markus-Epstein.
\newblock Stallings foldings and subgroups of amalgams of finite groups.
\newblock {\em Internat. J. Algebra Comput.}, 17(8):1493--1535, 2007.

\bibitem{2014:Nicaud}
C.~Nicaud.
\newblock Random deterministic automata.
\newblock In E.~Csuhaj{-}Varj{\'{u}}, M.~Dietzfelbinger, and Z.~{\'{E}}sik,
  editors, {\em Mathematical Foundations of Computer Science 2014 --- 39th
  International Symposium, {MFCS} 2014, Budapest, Hungary, August 25-29, 2014.
  Proceedings, Part {I}}, volume 8634 of {\em Lecture Notes in Computer
  Science}, pages 5--23. Springer, 2014.

\bibitem{1991:Short}
H.~Short.
\newblock Quasiconvexity and a theorem of {H}owson's.
\newblock In {\em Group theory from a geometrical viewpoint ({T}rieste, 1990)},
  pages 168--176. World Sci. Publ., River Edge, NJ, 1991.

\bibitem{1983:Stallings}
J.~R. Stallings.
\newblock Topology of finite graphs.
\newblock {\em Invent. Math.}, 71(3):551--565, 1983.

\bibitem{2006:Touikan}
N.~W.~M. Touikan.
\newblock A fast algorithm for {S}tallings' folding process.
\newblock {\em Internat. J. Algebra Comput.}, 16(6):1031--1045, 2006.

\end{thebibliography}
}

\end{document}